\documentclass[a4paper,dvipsnames,fleqn]{article}

\usepackage{amsmath}
\usepackage{amssymb}
\usepackage{amsthm}
\usepackage{array}
\usepackage[backend=biber,bibencoding=utf8,bibstyle=alphabetic,citestyle=alphabetic,sorting=nyt,maxnames=4]{biblatex}
	\AtEveryBibitem{\clearfield{doi}}
	\AtEveryBibitem{\clearfield{isbn}}
	\AtEveryBibitem{\clearfield{issn}}
	\AtEveryBibitem{\clearfield{url}}
	
	\renewbibmacro{in:}{\ifentrytype{article}{}{\printtext{\bibstring{in}\intitlepunct}}}
	\AtEveryBibitem{\ifentrytype{book}{\clearfield{pages}}{}}
\usepackage[colorlinks=true,allcolors=purple]{hyperref}

\usepackage{cleveref}
	\crefformat{equation}{#2(#1)#3}
\usepackage{dsfont}
\usepackage{enumitem}
\usepackage[OT2,OT1]{fontenc}
\usepackage[cal=cm,scr=boondoxo]{mathalfa}
\usepackage{pifont}
\usepackage{appendix}

\usepackage[textwidth=3cm,colorinlistoftodos]{todonotes}

\usepackage{xfrac}

\numberwithin{equation}{section}			% Enumeration of equations
\swapnumbers						% Enumeration of theorem-like environments

\newcommand\cyr{%					% cyrillic font
\renewcommand\rmdefault{wncyr}%
\renewcommand\sfdefault{wncyss}%
\renewcommand\encodingdefault{OT2}%
\normalfont
\selectfont}
\DeclareTextFontCommand{\textcyr}{\cyr}

\newcounter{Enum}				% Enumerated list
\newenvironment{Enumerate}{\begin{enumerate}[label={\rm ({\roman*})}]}{\end{enumerate}}

\newcommand{\descriptionlabelsave}{}		% Itemized list

\newcounter{StepsCount}				% Enumerated list with no indentation (e.g. steps in proof)

\newcounter{StepsRefCount}

\theoremstyle{plain}
	\newtheorem{lemma}{Lemma}[section]
	\newtheorem{proposition}[lemma]{Proposition}
	\newtheorem{theorem}[lemma]{Theorem}
	\newtheorem{subtheorem}{Theorem}[lemma]

	\newtheorem{corollary}[lemma]{Corollary}
	\newcommand{\GenericTheoremName}{}\newtheorem{generictheorem}[lemma]{\GenericTheoremName}
\theoremstyle{definition}
	\newtheorem{definition}[lemma]{Definition}
	\newcommand{\GenericDefinitionName}{}\newtheorem{genericdefinition}[lemma]{\GenericDefinitionName}
\theoremstyle{remark}
	\newtheorem{remark}[lemma]{Remark}
	\newtheorem{example}[lemma]{Example}

	\newcommand{\GenericRemarkName}{}\newtheorem{genericremark}[lemma]{\GenericRemarkName}
\newcommand{\mc}[1]{{\mathcal{#1}}}			% --- abbreviation ---
\newcommand{\ms}[1]{{\mathscr{#1}}}			% --- abbreviation ---
\newcommand{\bb}[1]{{\mathbb{#1}}}			% --- abbreviation ---
\newcommand{\mr}{\mathring}				% --- abbreviation ---
\newcommand{\Poi}[1]{\mc P[{#1}]}

\DeclareMathOperator{\RE}{Re}				% real part

\DeclareMathOperator{\IM}{Im}				% imaginary part

\DeclareMathOperator{\tr}{tr}	
\DeclareMathOperator{\sgn}{sgn}

\newcommand{\Side}[1]{\hfill{#1}\kern10pt}		% text put on the right side of line with offset

\newcommand{\smmatrix}[4]{\Bigl(			% small matrix for use in textline
\begin{smallmatrix}
\hspace*{-0.2ex} #1 \hspace*{0.2ex} & \hspace*{0.2ex} #2 \hspace*{-0.2ex}
\\[0.5ex]
\hspace*{-0.2ex} #3 \hspace*{0.2ex} & \hspace*{0.2ex} #4 \hspace*{-0.2ex}
\end{smallmatrix}
\Bigr)}
\newcommand{\Dummy}{\text{\textvisiblespace\kern1pt}}	% Platzhaltersymbol fuer Funktionsargumente
\newcommand{\Smallo}{{\rm o}}				% small o
\newcommand{\BigO}{{\rm O}}				% big o

\newcommand{\DS}{\mid\mkern3mu}				% delimiter for set definition
\newcommand{\DP}{{.\kern5pt}}				% delimiter for predicate formula
\newcommand{\DE}{\mathrel{\mathop:}=}			% defining equality
\newcommand{\DD}{\mkern4mu\mathrm{d}}			% distance and rm-d for integration differential

\begin{document}

\begin{flushleft}
	{\Large\bf	A quantitative formula for the imaginary part of a Weyl coefficient\\[5mm]
	
%	Quantitative formulae for Weyl coefficients: accessing tails of the spectral measure
}
	%\\[5mm]
	\textsc{
	Jakob Reiffenstein
		\hspace*{-14pt}
		\renewcommand{\thefootnote}{\fnsymbol{footnote}}
		\setcounter{footnote}{2}
			\footnote{
	Department of Mathematics, University of Vienna \\
	Oskar-Morgenstern-Platz 1, 1090 Wien, AUSTRIA \\
	email: jakob.reiffenstein@univie.ac.at}
	} \\[1ex]

\end{flushleft}

%%%%%%%%%%%%%%%%%%%%%%%%%%%%%%%%%%%%%%%%%%%%%%%%%%%%%%%%%%%%%%%%%%%%%%%%%%%%%%%%%%%%%%%%%
% disable in final version
%\pagenumbering{roman}
%\fbox{
%\parbox{100mm}{
%\hspace*{0pt}\\[1mm]
%\centerline{{\Large\ding{45}}\quad\,{\large\sc Draft}\quad{\Large\ding{45}}}
%\hspace*{0pt}\\[-2mm]
%\textcircledP\ \ Preliminary version
%\\[2mm]
%\ding{233}\quad Use pdflatex/biber to compile
%\\[-1mm]
%}
%}
%\\[1ex]
	{\small
	\textbf{Abstract.}
We investigate two-dimensional canonical systems $y'=zJHy$ on an interval, with positive semi-definite Hamiltonian $H$%, such that limit circle case prevails at the left endpoint and limit point case at the right
. Let $q_H$ be the Weyl coefficient of the system. We prove a formula that determines the imaginary part of $q_H$ along the imaginary axis up to multiplicative constants, which are independent of $H$. We also provide versions of this result for Sturm-Liouville operators and Krein strings. \\
Using classical Abelian-Tauberian theorems, we deduce characterizations of spectral properties such as integrability of a given comparison function w.r.t. the spectral measure $\mu_H$, and boundedness of the distribution function of $\mu_H$ relative to a given comparison function. \\
We study in depth Hamiltonians for which $\arg q_H(ir)$ approaches $0$ or $\pi$ (at least on a subsequence). It turns out that this behavior of $q_H(ir)$ imposes a substantial restriction on the growth of $|q_H(ir)|$. Our results in this context are interesting also from a function theoretic point of view.
	\\[3mm]
	\textbf{AMS MSC 2020:} 30E99, 34B20, 34L05, 34L40
	\\
	\textbf{Keywords:} Canonical system, Weyl coefficient, growth estimates, high-energy behaviour
	}

\pagenumbering{arabic}
\setcounter{page}{1}
\setcounter{footnote}{0} 
%
%%%%%%%%%%%%%%%%%%%%%%%%%%%%%%%%%%%%%%%%%%%%%%%%%%%%%%%%%%%%%%%%%%%%%%%%%%%%%%%%%%%%%%%%%

%---------
%   TEXTBODY
%---------

%
%
%
%\section[\textcolor{ForestGreen}{}]{}
%\section[\textcolor{Dandelion}{}]{}
%\section[\textcolor{BrickRed}{...}]{...}
\section[{Introduction}]{Introduction}

\noindent We study two-dimensional \textit{canonical systems} 
\begin{align}
\label{A33}
y'(t)=zJH(t)y(t), \quad \quad t \in [a,b) \, \text{ a.e.},
\end{align}
where $-\infty < a <b \leq \infty$, $z \in \bb C$ is a spectral parameter and $J:=\smmatrix 0{-1}10$. The \textit{Hamiltonian} $H$ is assumed to be a locally integrable, $\bb R^{2 \times 2}$-valued function on $[a,b)$ that further satisfies
\begin{itemize}
\item[$\rhd$] $H(t) \geq 0$ and $H(t) \neq 0$, \quad \quad $t \in [a,b)$ a.e.;
\item[$\rhd$] $H$ is definite, i.e., if $v \in \mathbb{C}^2$ is s.t. $H(t)v \equiv 0$ on $[a,b)$, then $v=0$;
\item[$\rhd$] $\int_a^b \tr H(t) \DD t=\infty$ (limit point case at $b$).
\end{itemize}
Together with a boundary condition at $a$, the equation (\ref{A33}) becomes the eigenvalue equation of a self-adjoint (possibly multi-valued) operator $A_H$ in a Hilbert space $L^2(H)$ associated with $H$. Throughout this paper, we fix the boundary condition $(1,0)y(a)=0$, which is no loss of generality. \\ 
Many classical second-order differential operators such as Schr\"odinger and Sturm-Liouville operators, Krein strings, and Jacobi operators can be transformed to the form (\ref{A33}), see, e.g., \cite{remling:2018,teschl:2009,behrndt.hassi.snoo:2020,kaltenbaeck.winkler.woracek:2007,kac:1999}. Canonical systems thus form a unifying framework. \\
All of the above operators have in common that their spectral theory is centered around the Weyl coefficient $q$ of the operator (also referred to as Titchmarsh-Weyl $m$-function). This function is constructed by Weyl's nested disk method and is a Herglotz function, i.e., it is holomorphic on $\bb C \setminus \bb R$ and satisfies there $\frac{\IM q(z)}{\IM z} \geq 0$ as well as $q(\overline{z})=\overline{q(z)}$. It can thus be represented as
\begin{align}
\label{A17}
q(z)=\alpha + \beta z + \int_{\bb R} \bigg(\frac{1}{t-z}-\frac{t}{1+t^2} \bigg) \DD \mu(t), \quad \quad z \in \bb C \setminus \bb R
\end{align}
with $\alpha \in \bb R$, $\beta \geq 0$, and $\mu$ a positive Borel measure on $\bb R$ satisfying $\int_{\bb R} \frac{d\mu(t)}{1+t^2} <\infty$. The measure $\mu$ in the integral representation (\ref{A17}) of the Weyl coefficient is a spectral measure of the underlying operator model if $\beta =0$ (if $\beta > 0$, a one-dimensional component has to be added). The importance of canonical systems in this context lies in the Inverse Spectral Theorem of L. de Branges, stating that each Herglotz function $q$ is the Weyl coefficient of a unique (suitably normalized) canonical system. \\

\noindent Given a Hamiltonian $H$, we are ultimately interested in the description of properties of its spectral measure $\mu_H$ in terms of $H$. The correspondence between $H$ and $\mu_H$ can be best understood using the Weyl coefficient $q_H$, whose imaginary part $\IM q_H$ determines $\mu_H$ via the Stieltjes inversion formula. \\
In their recent paper \cite{langer.pruckner.woracek:heniest}, Langer, Pruckner, and Woracek gave a two-sided estimate for $\IM q_H(ir)$ in terms of the coefficients of $H$:
\begin{align}
\label{A43}
L(r) \lesssim \IM q_H(ir) \lesssim A(r), \quad \quad r>0,
\end{align}
where $L,A$ are explicit in terms of $H$, and we used the notation $f(r) \lesssim g(r)$ to state that $f(r) \leq Cg(r)$ for a constant $C>0$. Moreover, in (\ref{A43}) the constants implicit in $\lesssim$ are independent of $H$. The exact formulation of this result will be recalled in \Cref{Y98}. \\
It may happen that $L(r)=\Smallo (A(r))$, and $\IM q_H(ir)$ is not determined by (\ref{A43}). A toy example for this is the Hamiltonian
\begin{align*}
H(t)=t \left(\begin{matrix}
|\log t|^{\color{white} 1} & |\log t|^2 \\
|\log t|^2 & |\log t|^3 \\
\end{matrix} \right), \quad \quad t \in [0,\infty).
\end{align*}
For $r \to \infty$, a calculation shows that
\begin{align*}
 L(r) &\asymp (\log r)^{-3}, \quad \quad A(r) \asymp (\log r)^{-1},
\end{align*}
where $f(r) \asymp g(r)$ means that both $f(r) \lesssim g(r)$ and $g(r) \lesssim f(r)$.
\newline

\noindent The following theorem, which is our main result, improves the estimate (\ref{A43}) by giving a formula for $\IM q_H(ir)$ up to universal multiplicative constants.

\begin{theorem}
\label{T1}
Let $H$ be a Hamiltonian on $[a,b)$, and denote\footnote{When there is no risk of ambiguity, we write $\Omega$ and $\omega_j$ instead of $\Omega_H$ and $\omega_j^{(H)}$ for short.}
\begin{equation}\label{Y08}
		H(t) = \begin{pmatrix} h_1(t) & h_3(t) \\ h_3(t) & h_2(t) \end{pmatrix},\quad
		\Omega_H(t) = \begin{pmatrix} \omega_1^{(H)}(t) & \omega_3^{(H)}(t) \\ \omega_3^{(H)}(t) & \omega_2^{(H)}(t) \end{pmatrix}
		\DE\int_a^t H(s)\DD s.
	\end{equation}
Let $\hat t : (0,\infty) \to (a,b)$ be a function satisfying\footnote{We will see later that the equation $\det \Omega_H(t)=\frac{1}{r^2}$ has a unique solution for every $r>0$. A possible choice of $\hat t$ is thus the function that maps $r>0$ to this solution. }
\begin{align} 
\label{A49}
\det \Omega_H(\hat t(r)) \asymp \frac{1}{r^2}, \quad \quad r \in (0,\infty).
\end{align}  
Then
\begin{align}
\label{A2}
\IM q_H(ir) &\asymp \bigg|q_H(ir)-\frac{\omega_3^{(H)}(\hat t(r))}{\omega_2^{(H)}(\hat t(r))} \bigg| \asymp \frac{1}{r\omega_2^{(H)}(\hat t(r))}, \\[1.7ex]
\label{A3}
\frac{\IM q_H(ir)}{|q_H(ir)|^2} &\asymp \frac{1}{r\omega_1^{(H)}(\hat t(r))},
\end{align} 
for $r \in (0,\infty)$. The constants implicit in $\asymp$ in (\ref{A2}) and (\ref{A3}) depend on the constants hidden in $\asymp$ in (\ref{A49}), but not on $H$. \\
If, in addition, $\IM q_H(ir)=\Smallo (|q_H(ir)|)$ for $r \to \infty$ (or $r \to 0$), then\footnote{With $f(r) \sim g(r)$ meaning $\lim \frac{f(r)}{g(r)}=1.$}
\begin{align}
\label{A11}
q_H(ir) \sim \frac{\omega_3^{(H)}(\hat t(r))}{\omega_2^{(H)}(\hat t(r))}, \quad \quad r \to \infty \quad ( r \to 0).
\end{align}
\end{theorem}

\noindent The two-sided estimate (\ref{A2}) has some useful features: its pointwise nature, its applicability for $r \to \infty$ and $r \to 0$, and the universality of the constants hidden in $\asymp$. However, it is rather different from an asymptotic formula: it does not capture small oscillations of $\IM q_H(ir)$ around $\frac{1}{r\omega_2^{(H)}(\hat t(r))}$. \\
Note also that the first relation in (\ref{A2}) can be seen as a statement about the real part of $q_H(ir)$. In fact, $\IM q_H(ir)$ is also obtained if we subtract $\RE q_H(ir)$ from $q_H(ir)$, then take absolute values. It is an open question whether $\RE q_H(ir)$ can be described more directly in terms of $H$.
\newline

\noindent A most important class of operators is that of Sturm-Liouville (in particular, Schr\"odinger) operators. Let us provide a reformulation of \Cref{T1} for these operators right away.

%the numbers $L(r)$ and $A(r)$ from (\ref{A41}) describe the imaginary parts of bottom and top of suitably parametrised Weyl disks containing $q_H(re^{i\vartheta})$. We will see 

 %related to the position of $q_H(ir)$ within Weyl disks. This is discussed in \Cref{S4}.
 
\subsection*{Sturm-Liouville operators}

We provide a version of \Cref{T1} for Sturm-Liouville equations
\begin{align}
\label{A44}
-(py')'+qy=zwy 
\end{align}
on $(a,b)$, where $1/p, q,w \in L^1_{loc}(a,b)$, $w>0$ and $p,q$ are real-valued. Suppose that $a$ is in limit circle case and $b$ is in limit point case. Impose a Dirichlet boundary condition at $a$, i.e., $y(a)=0$. The Weyl coefficient for this problem is the unique number $m(z)$ with
\[
c(z,\cdot)+m(z)s(z,\cdot) \in L^2((a,b),w(x)\DD x)
\]
where $c(z,\cdot)$ and $s(z,\cdot)$ are solutions of (\ref{A44}) with initial values 
\[
\binom{p(a)c'(z,a)}{c(z,a)}=\binom{0}{1}, \quad \binom{p(a)s'(z,a)}{s(z,a)}=\binom{1}{0}.
\]
\begin{theorem}
\label{T9}
For each $t \in (a,b)$, let $(.,.)_t$ and $\|.\|_t$ denote the scalar product and norm on $L^2((a,t),w(x)\DD x)$, i.e.,
\[
(f,g)_t=\int_a^t f(x)\overline{g(x)} w(x) \DD x.
\]
For $\xi \in \mathbb{R}$, let $\hat t_\xi : (0,\infty) \to (a,b)$ be a function satisfying
\begin{align}
\label{A51}
\|c(\xi,\cdot)\|_{\hat t_\xi(r)}^2 \|s(\xi,\cdot)\|_{\hat t_\xi(r)}^2 - (c(\xi,\cdot),s(\xi,\cdot))_{\hat t_\xi(r)}^2 \asymp \frac{1}{r^2}, \quad \,\, r \in (0,\infty).
\end{align}
Then
%\bigg|m(\xi+ir)+\frac{(c(\xi,\cdot),s(\xi,\cdot))_{\hat t_\xi(r)}}{\|s(\xi,\cdot)\|_{\hat t_\xi(r)}^2} \bigg| \asymp
\begin{align}
\label{A45}
\IM m(\xi+ir) &\asymp \frac{1}{r \|s(\xi,\cdot)\|_{\hat t_\xi(r)}^2}, \\
\label{A46}
\frac{\IM m(\xi+ir)}{|m(\xi+ir)|^2} &\asymp \frac{1}{r \|c(\xi,\cdot)\|_{\hat t_\xi(r)}^2},
\end{align}
for $r \in (0,\infty)$. The constants implicit in $\asymp$ are independent of $p,q,w$ as well as $\xi$, but do depend on the constants pertaining to $\asymp$ in (\ref{A51}).
\end{theorem}

\noindent In fact, \Cref{T9} is a direct consequence of \Cref{T1} upon employing a transformation (cf. \cite{remling:2018} for $p=w=1$ and $\xi=0$) that maps solutions of (\ref{A44}) to solutions of the canonical system $y'=(z-\xi)JH_\xi y$, where 
\[
H_\xi (t) = w(t) \cdot \begin{pmatrix}
c(\xi,t)^2 & -s(\xi,t)c(\xi,t) \\
-s(\xi,t)c(\xi,t) & s(\xi,t)^2
\end{pmatrix}, \quad \quad t \in [a,b).
\]
The Weyl coefficients then satisfy $m(z)=q_{H_\xi}(z-\xi)$.

\subsection*{Historical remarks}

\noindent The origins of the Weyl coefficient in the theory of the Sturm-Liouville differential equation are well summarized in Everitt's paper \cite{everitt:2004}. We give a short account specifically on the history of estimates for the growth of the Weyl coefficient, which date back at least to the 1950s. Particular attention was often given to the deduction of asymptotic formulae for the Weyl coefficient \cite{marchenko:1952,kac:1973a,everitt:1972,kasahara:1975,atkinson:1981,bennewitz:1989}. However, asymptotic results usually depend on rather strong assumptions on the data. When weakening these assumptions, one can still ask for explicit estimates for $q(z)$ as $z \to \infty$ nontangentially in the upper half-plane. There is a number of rather early results that determine $|q(z)|$ up to $\asymp$, e.g., \cite{hille:1963,atkinson:1988,bennewitz:1989}, although these still depend on data subject to additional restrictions. Fundamental progress has been made by Jitomirskaya and Last \cite{jitomirskaya.last:1999}, who considered Schr\"odinger operators with arbitrary (real-valued and locally integrable) potentials. They found a formula up to $\asymp$ for $|q(z)|$, which also covers the case $z \to 0$. An analog of this formula for canonical systems was given in \cite{hassi.remling.snoo:2000}. \\
When it comes to $\IM q(z)$, however, no such formula was available. Only the very recent estimate (\ref{A43}) from \cite[Theorem 1.1]{langer.pruckner.woracek:heniest} made it possible to obtain our main result that determines $\IM q(z)$ up to $\asymp$.

\subsection*{Structure of the paper}

\noindent The proof of \Cref{T1}, together with some immediate corollaries, makes up \Cref{S2}. In \Cref{S5}, we continue with a first application, a criterion for integrability of a given comparison function with respect to $\mu_H$. We also characterize boundedness of the distribution function of $\mu_H$ relative to a given comparison function.
\newline

\noindent \Cref{S3} is dedicated to the boundary behavior of Herglotz functions. Cauchy integrals and the relative behavior of its imaginary and real part have been intensively studied. For example, for a Herglotz function $q$ it is known \cite{poltoratski:2003} that the set of $\xi \in \mathbb{R}$ for which 
\begin{align}
\label{A47}
\lim_{r \to 0} \frac{\IM q(\xi+ir)}{|q(\xi+ir)|}=0
\end{align}
is a zero set w.r.t. $\mu$. In contrast to measure theoretic results like this, we use the de Branges correspondence $H \leftrightarrow q_H$ to investigate this behavior pointwise w.r.t. $\xi$. In \Cref{T7} we show that if $\xi$ is such that (\ref{A47}) holds, then $|q(\xi+ir)|$ is slowly varying (cf. \Cref{A48}). \Cref{T8} is a partial converse of this statement.
\newline

\noindent In \Cref{S4} we turn to a finer study of $\IM q_H(ir)$ in the context of the geometric origins of (\ref{A43}) and (\ref{A2}). Namely, the functions $L$ and $A$ describe the imaginary parts of bottom and top of certain Weyl disks containing $q_H(ir)$. We show that there are restrictions on the possible location of $q_H(ir)$ within the disks, and construct a Hamiltonian $H$ for which $q_H(ir)$ oscillates back and forth between the bottoms and tops of the disks. This construction allows us to answer several open problems that were posed in \cite{langer.pruckner.woracek:heniest}.
\newline

\noindent We conclude our work with a reformulation of \Cref{T1} for the principal Titchmarsh-Weyl coefficient $q_S$ of a Krein string. This reformulation is the content of \Cref{S6}. 
%The reformulation yields a formula up to multiplicative constants for $\IM q_S(z)$ along the imaginary axis. Let us point out that for $z \in (-\infty ,0)$, $q_S(z)$ was already determined up to multiplicative constants in \cite{langer.pruckner.woracek:heniest} and, for well-behaved $\mathfrak{m}$, in \cite{kasahara.watanabe:2010}. Results on the spectral measure of a string, similar to what we do for canonical systems in \Cref{S5}, can already be found in \cite{langer.pruckner.woracek:heniest}.

\subsection*{Notation associated to Hamiltonians}

\noindent Let $H$ be a Hamiltonian on $[a,b)$. 
\newline

\noindent An interval $(c,d) \subseteq [a,b)$ is called $H$-\textit{indivisible} if $H(t)$ takes the form $h(t)\binom{\cos \varphi}{\sin \varphi}\binom{\cos \varphi}{\sin \varphi}^*$ a.e. on $(c,d)$, with scalar-valued $h$ and fixed $\varphi \in [0,\pi)$. The angle $\varphi$ is then called the \textit{type} of the interval.

\begin{definition}
Let \begin{align}
\mr a (H) &:=\inf \Big\{t > a \,\Big| \, (a,t) \text{ is not }H\text{-indivisible of type } 0 \text{ or } \frac{\pi}{2} \Big\}, \\
\hat a (H) &:=\inf \Big\{t > a \,\Big|\, (a,t) \text{ is not }H\text{-indivisible} \Big\}.
\end{align}
Usually, we write $\mr a$ and $\hat a$ for short. Since $H$ is assumed to be definite, both of these numbers are smaller than $b$.
\end{definition}
\noindent Note that $(\omega_1 \omega_2)(t)>0$ if and only if $(a,t)$ is not $H$-indivisible of type $0$ or $\frac{\pi}{2}$, i.e., $t>\mr a$. Using the assumption $\int_a^b \tr H(t) \DD t=\infty$, we infer that $\omega_1 \omega_2$ is an increasing bijection from $(\mr a,b)$ to $(0,\infty)$. \\
Similarly, $\det \Omega (t)>0$ is equivalent to $t>\hat a$. We have
\[
\frac{d}{dt} \Big(\frac{\det \Omega (t)}{\omega_1(t)} \Big)=\omega_1(t)^{-2} \binom{-\omega_3(t)}{\omega_1(t)}^* H(t)\binom{-\omega_3(t)}{\omega_1(t)} \geq 0
\]
and (by symmetry) $\frac{d}{dt} \big(\frac{\det \Omega}{\omega_2}\big) \geq 0$. Since at least one of $\omega_1$ and $\omega_2$ is unbounded, $\det \Omega$ is an increasing bijection from $(\hat a,b)$ to $(0,\infty)$.

\begin{definition}
\label{A35}
For a Hamiltonian $H$ and a number $\eta >0$, set \\
\begin{minipage}{.5\linewidth}
\begin{equation*}
\mr r_{\eta,H} : \left\{\begin{array}{ccc}
(\mr a,b) &\to &(0,\infty) \\[0.5ex]
t &\mapsto & \frac{\eta}{ 2\sqrt{(\omega_1 \omega_2)(t)}},
\end{array}\right.
\end{equation*}
\end{minipage}%
\begin{minipage}{.5\linewidth}
\begin{equation*}
\hat r_{\eta,H} : \left\{\begin{array}{ccc}
(\hat a,b) &\to &(0,\infty) \\[0.5ex]
t &\mapsto & \frac{\eta}{ 2\sqrt{\det \Omega (t)}}.
\end{array}\right.
\end{equation*}
\end{minipage}
\vspace{1ex}

\noindent Both of these functions are decreasing and bijective. We define their inverse functions,
\begin{align}
\mr t_{\eta,H}:=\mr r_{\eta,H}^{-1} \,:\, (0,\infty) \to (\mr a,b), \quad \quad \hat t_{\eta,H}:=\hat r_{\eta,H}^{-1} \,:\, (0,\infty) \to (\hat a,b).
\end{align}
Note that the functions $\hat t_{\eta,H}$, for any $\eta>0$, satisfy (\ref{A49}). Functions of this form will be the default choice of $\hat t$ for the sake of \Cref{T1}. We will often fix $\eta$ and $H$ and write $\mr r$, $\mr t$, $\hat r$, $\hat t$ for short. If $\eta$ is fixed but the Hamiltonian is ambiguous, we may write $\mr r_H$, $\mr t_H$, $\hat r_H$, $\hat t_H$ to indicate dependence on $H$.
\end{definition}

%\noindent In other words, for $r>0$ and $\eta$ fixed, $\mr t_{\eta,H}(r)$ and $\hat t_{\eta,H}(r)$ are the unique numbers satisfying 
%\begin{align}
%(\omega_1 \omega_2)(\mr t_{\eta,H}(r))=\frac{\eta ^2}{4r^2}, \quad \det \Omega (\hat t_{\eta,H}(r))=\frac{\eta ^2}{4r^2},
%\end{align}
%respectively.

\section{On the imaginary part of the Weyl coefficient}
\label{S2}

\noindent We start by providing the details of the estimate (\ref{A43}), which is the central result in \cite{langer.pruckner.woracek:heniest}.

\begin{theorem}[{\cite[Theorem 1.1]{langer.pruckner.woracek:heniest}}]
\label{Y98}
Let $H$ be a Hamiltonian on $[a,b)$, and let $\eta \in (0,1-\frac{1}{\sqrt 2})$ be fixed. For $r>0$, let $\mr t(r)$ be the unique number satisfying
\begin{align}
\label{A50}
(\omega_1^{(H)} \omega_2^{(H)})(\mr t(r))=\frac{\eta^2}{4r^2},
\end{align} 
cf. \Cref{A35}. Set\footnote{If $\eta$ and $H$ are clear from the context, we may write $A$ and $L$ for short.}
\[
A_{\eta,H}(r):=\frac{\eta}{2r\omega_2^{(H)}(\mr t(r))}, \quad \quad L_{\eta,H}(r):= \frac{\det \Omega_H (\mr t(r))}{(\omega_1^{(H)} \omega_2^{(H)})(\mr t(r))} \cdot A_{\eta,H}(r).
\]
Then the Weyl coefficient $q_H$ associated with the Hamiltonian $H$ satisfies
\begin{align}
|q_H(ir)| &\asymp A_{\eta,H}(r),
\label{Y35} \\[1.5ex]
L_{\eta,H}(r) \lesssim \IM q_H(ir) &\lesssim A_{\eta,H}(r)
\label{Y96}
\end{align}
for $r \in (0,\infty)$. The constants implicit in these relations are independent of $H$. Their dependence on $\eta$ is continuous.
\end{theorem}
\noindent In the following proof of \Cref{T1}, we will also show that \Cref{Y98} still holds if $\mr t: (0,\infty) \to (a,b)$ is a function satisfying $(\omega_1 \omega_2)(\mr t(r)) \asymp \frac{1}{r^2}$, and
\[
A(r):=\frac{1}{r\omega_2^{(H)}(\mr t(r))}, \quad \quad L(r):= \frac{\det \Omega_H (\mr t(r))}{(\omega_1^{(H)} \omega_2^{(H)})(\mr t(r))} \cdot A(r).
\]
In particular, we can choose any $\eta>0$ in (\ref{A50}).

\begin{proof}[Proof of \Cref{T1}]

Let $\hat t_{\eta,H}$ be defined as in \Cref{A35}. We show that for any $\eta>0$, \Cref{T1} holds for $\hat t_{\eta,H}$ in place of $\hat t$, and that the dependence on $\eta$ of the constants hidden in $\asymp$ in (\ref{A2}) and (\ref{A3}) is continuous. This then implies that \Cref{T1} holds for any function $\hat t$ satisfying (\ref{A49}).
\newline

\noindent The proof is divided into steps. \\
\item[\textbf{Step 1.}] 
We introduce a family of transformations of $H$ that leave the imaginary part of the Weyl coefficient unchanged. If $p \in \bb R$ and 
\[
H_p(t):=\smmatrix 1p01 H(t) \smmatrix 10p1 = 
\begin{pmatrix}
h_1(t)+2p h_3(t)+p^2 h_2(t)  & h_3(t)+ph_2(t) \\
h_3(t)+ph_2(t) & h_2(t)
\end{pmatrix},
\] 
an easy calculation shows that the Weyl coefficient $q_p$ of $H_p$ is given by $q_p(z)=q_0(z)+p=q_H(z)+p$. \\

\item[\textbf{Step 2.}] 
We prove (\ref{A2})-(\ref{A11}) for fixed $\eta \in (0,1-\frac{1}{\sqrt{2}})$. The following abbreviations are used only in Step 2: \\

\begin{tabular}{|r|l||r|l||r|l|@{}m{0pt}@{}}
\hline
short form & meaning & short form & meaning &short form & meaning & \\[10pt]
\hline \hline 
\rule{0pt}{3ex}$\mr t$  & \rule{0pt}{3ex} $\mr t_{\eta,H}$ & \rule{0pt}{3ex} $\mr t_p$ & \rule{0pt}{3ex} $\mr t_{\eta,H_p}$ & \rule{0pt}{3ex} $\Omega_p$ & \rule{0pt}{3ex} $\Omega_{H_p}$&  \\[3pt]
\hline
\rule{0pt}{3ex}$\hat t$ & \rule{0pt}{3ex}$\hat t_{\eta,H}$ & \rule{0pt}{3ex}$\hat t_p$ & \rule{0pt}{3ex}$\hat t_{\eta,H_p}$ & \rule{0pt}{3ex}$\omega_j^{(p)}$ & \rule{0pt}{3ex}$\omega_j^{(H_p)}$ &  \\[3pt]
\hline 
\rule{0pt}{3ex}$L_p$& \rule{0pt}{3ex}$L_{\eta,H_p}$ & \rule{0pt}{3ex}$A_p$& \rule{0pt}{3ex}$A_{\eta,H_p}$ & \rule{0pt}{3ex}$\Omega$ & \rule{0pt}{3ex}$\Omega_H$ 
& \\[3pt]
\hline
\end{tabular}
\\[4pt]

%\begin{itemize}
%\item[$\rhd$] $\mr t$ and $\hat t$ are short for $\mr t_{\eta,H}$ and $\hat t_{\eta,H}$, cf. \Cref{A35}.
%\item[$\rhd$] $\mr t_p$ and $\hat t_p$ are short for $\mr t_{\eta,H_p}$ and $\hat t_{\eta,H_p}$.
%\item[$\rhd$] $\Omega$ and $\omega_j$ abbreviate $\Omega_H$ and $\omega_j^{(H)}$.
%\item[$\rhd$] $\Omega_p$ and $\omega_j^{(p)}$ stand for $\Omega_{H_p}$ and $\omega_j^{(H_p)}$.
%\item[$\rhd$] $L_p$ and $A_p$ are as in \Cref{Y98}, but with $H_p$ in place of $H$.
%\end{itemize}
\noindent Let $r>0$ be fixed (this is important). Our first observation is that $\hat t_p(r)=\hat t(r)$ for any $p$ since $\det \Omega_p(t)=\det \Omega (t)$ does not depend on $p$. If we can find $p$ such that $\mr t_p(r)=\hat t_p(r)=\hat t(r)$, then clearly 
\[
\frac{L_p(r)}{A_p(r)}=\frac{\det \Omega_p(\mr t_p(r))}{(\omega_1^{(p)}\omega_2^{(p)})(\mr t_p(r))}=\frac{\det \Omega_p(\hat t_p(r))}{(\omega_1^{(p)}\omega_2^{(p)})(\mr t_p(r))}=1.
\]
We apply \Cref{Y98} with $\eta$ and $H_p$. The estimate (\ref{Y96}) then takes the form
\begin{equation}
\label{P1}
A_p(r) = L_p(r) \lesssim \IM q_H(ir) \lesssim A_p(r)
\end{equation}
while (\ref{Y35}) turns into
\begin{equation}
\label{P2}
 |q_H(ir)+p| \asymp A_p(r),
\end{equation}
where
\[
A_p(r)= \frac{\eta}{2r\omega_2^{(p)}(\mr t_p(r))} = \frac{\eta}{2r\omega_2(\hat t(r))}.
\]
The right choice of $p$ is
\[
p=- \frac{\omega_3(\hat t(r))}{\omega_2(\hat t(r))},
\]
leading to $\omega_3^{(p)}(\hat t(r))=0$ and thus
\[
(\omega_1^{(p)}\omega_2^{(p)})(\hat t(r))=\det \Omega_p (\hat t(r))=\det \Omega (\hat t(r))=\frac{\eta^2}{4r^2}.
\]
Consequently, $\mr t_p(r)=\hat t(r)$. Observe that the implicit constants in (\ref{Y35}) and (\ref{Y96}) are independent of $H$ and $r$ and depend continuously on $\eta$. This shows that (\ref{A2}) holds, with constants depending continuously on $\eta$. \\
\item[\textbf{Step 3.}] (\ref{A3}) follows from an application of (\ref{A2}) to $\tilde H:=J^{\top}HJ=\smmatrix {h_2}{-h_3}{-h_3}{h_1}$ and note that $\hat t_{\eta,\tilde H}=\hat t_{\eta,H}$. Thus
\[
\frac{\IM q_H(ir)}{|q_H(ir)|^2}=\IM \Big(- \frac{1}{q_H(ir)} \Big)= \IM q_{\tilde H}(ir) \asymp \frac{1}{r \omega_2^{(\tilde H)}(\hat t_{\eta,\tilde H}(r))}=\frac{1}{r \omega_1^{(H)}(\hat t_{\eta,H}(r))}.
\]
Formula (\ref{A11}) follows if we divide (\ref{A2}) by $|q_H(ir)|$. Hence, we proved the assertion for $\eta \in (0,1-\frac{1}{\sqrt{2}})$.
\newline

\noindent In the remaining steps we treat the missing case $\eta \geq 1-\frac{1}{\sqrt{2}}$.
\item[\textbf{Step 4.}] Let $k>0$. For use in Step 5, we show that
\begin{align}
\label{A31}
\IM q_H(ir) \asymp \IM q_H (ikr), \quad \quad |q_H(ir)| \asymp |q_H (ikr )|
\end{align}
for $r \in (0,\infty)$, where the constants in $\asymp$ depend continuously on $k$ and are independent of $H$. \\
For the imaginary part, the statement is easy to see from the integral representation (\ref{A17}). For the absolute value, we use the Hamiltonian $\tilde H$ from Step 3 to obtain
\[
\frac{\IM q_H(ir)}{|q_H(ir)|^2}=\IM q_{\tilde H}(ir) \asymp \IM q_{\tilde H}(ikr)=\frac{\IM q_H(ikr)}{|q_H(ikr)|^2}.
\]
This shows that $|q_H(ir)| \asymp |q_H (ikr )|$ as well.
%To prove (\ref{A31}), we use the notation $\mr t_{\eta,H}(r)$, $\hat t_{\eta,H}(r)$, to indicate dependence on $\mr t(r)$ and $\hat t(r)$ on the choice of parameter $\eta \in (0,1-\frac{1}{\sqrt{2}})$ and Hamiltonian $H$. We choose $\eta >0$ such that both $\eta$ and $k\eta $ are in $(0,1-\frac{1}{\sqrt{2}})$. Upon noting that
%\[
%\mr t_{\eta,H}(r)=\mr t_{k\eta,kH}(r), \quad \quad \hat t_{\eta,H}(r)=\hat t_{k\eta,kH}(r)
%\]
%for all $r>0$, we obtain from \Cref{Y98} that
%\[
%|q_H(ir)| \asymp \frac{1}{r\omega_2(\mr t_{\eta,H}(r))} =k \cdot \frac{1}{r(k\omega_2)(\mr t_{k\eta,kH}(r))} \asymp |q_H(ikr)|
%\]
%and from (\ref{A2}) that
%\[
%\IM q_H(ir) \asymp \frac{1}{r\omega_2(\hat t_{\eta,H}(r))} =k \cdot \frac{1}{r(k\omega_2)(\hat t_{k\eta,kH}(r))} \asymp \IM q_H(ikr).
%\]

\item[\textbf{Step 5.}]
\noindent Fix a Hamiltonian $H$, and let $\eta_0 \geq 1-\frac{1}{\sqrt{2}}$. Then
\[
\mr t_{\eta_0,H}(r)=\mr t_{\frac 14,\frac{1}{4\eta_0}H}(r), \quad \quad \hat t_{\eta_0,H}(r)=\hat t_{\frac 14,\frac{1}{4\eta_0}H}(r)
\]
and
\[
A_{\eta_0,H}(r)= A_{\frac 14,\frac{1}{4\eta_0}H}(r), \quad \quad L_{\eta_0,H}(r)= L_{\frac 14,\frac{1}{4\eta_0}H}(r).
\]

Since $\frac 14$ is less than $1-\frac{1}{\sqrt{2}}$, we can use \Cref{Y98} with $\eta:=\frac 14$ to obtain
\begin{align}
\label{A30}
L_{\eta_0,H}(r) = L_{\frac 14,\frac{1}{4\eta_0}H}(r) \lesssim &\IM q_{\frac{1}{4\eta_0}H} (ir) \\
\leq &|q_{\frac{1}{4\eta_0}H}(ir)| \asymp A_{\frac 14,\frac{1}{4\eta_0}H}(r) = A_{\eta_0,H}(r) \nonumber
\end{align}
for $r \in (0,\infty)$. Since $q_{\frac{1}{4\eta_0}H}(z)=q_H \big(\frac{z}{4\eta_0} \big)$ and by Step 4, we see that \Cref{Y98} holds for arbitrary $\eta >0$. It is easy to check that continuous dependence of constants on $\eta$ is retained. Repeating Steps $1-3$ now shows that also \Cref{T1} holds for $\hat t_{\eta,H}$ for any $\eta >0$. Moreover, it is not hard to see that everything still works if $\hat t$ is a function satisfying (\ref{A49}).
\end{proof}

%\begin{remark}
%In Step 4 of the previous proof, we obtained the following statement: \\
%\textit{Let $q$ be a Herglotz function, and let $k>0$. Then} 
%\begin{align*}
%|q(ir)| \asymp |q (ikr )|, \quad \quad \IM q(ir) \asymp \IM q(ikr)
%\end{align*}
%\textit{for $r \in (0,\infty)$, with constants independent of $q$.} \\
%The change in perspective, from Hamiltonians to Herglotz functions, is possible due to the fact that each Herglotz function $q$ is the Weyl coefficient of some Hamiltonian $H$. 
%\end{remark}

\begin{remark}
\label{A36}
\Cref{Y98} and \Cref{T1}, in the form we stated them, give information about $q_H(z)$ for $z=ir$. However, if $\vartheta \in (0,\pi )$ is fixed, these theorems also hold
\begin{itemize}
\item [$\rhd$] for $z=re^{i\vartheta}$ uniformly for $r \in (0,\infty )$ and
\item [$\rhd$] for $z=re^{i\varphi}$ uniformly for $r \in (0,\infty)$ and $|\frac{\pi}{2}-\varphi | \leq |\frac{\pi}{2}-\vartheta |$.
\end{itemize}
We restate the explicit constants coming from \cite{langer.pruckner.woracek:heniest}. Fix $\eta \in (0,1-\frac{1}{\sqrt{2}})$ and set $\sigma :=(1-\eta )^{-2}-1 \in (0,1)$. With
\begin{align*}
c_-(\eta ,\vartheta)=\frac{\eta \sin \vartheta}{2(1+|\cos \vartheta |)} \cdot \frac{1-\sigma}{1+\sigma}, \quad \quad c_+(\eta ,\vartheta)=\frac{\sigma+\frac{2}{\eta \sin \vartheta}}{1-\sigma},
\end{align*}
we have\footnote{Since $c_-$ and $c_+$ are clearly monotonic in $\vartheta $, (\ref{A37}) and (\ref{A38}) still hold when $q_H(re^{i\vartheta})$ is replaced by $q_H(re^{i\varphi})$, where $|\frac{\pi}{2}-\varphi | \leq |\frac{\pi}{2}-\vartheta |$.}
\begin{align}
\label{A37}
c_-(\eta ,\vartheta) \cdot \frac{\eta}{2} \cdot \frac{1}{r\omega_2(\hat t_{\eta ,H}(r))} &\leq \IM q_H(re^{i\vartheta}) \leq c_+(\eta ,\vartheta) \cdot \frac{\eta}{2} \cdot \frac{1}{r\omega_2(\hat t_{\eta ,H}(r))}, \\[1ex]
\label{A38}
c_-(\eta ,\vartheta) \cdot \frac{\eta}{2} \cdot \frac{1}{r\omega_1(\hat t_{\eta ,H}(r))} &\leq \frac{\IM q_H(re^{i\vartheta})}{|q_H(re^{i\vartheta})|^2} \leq c_+(\eta ,\vartheta) \cdot \frac{\eta}{2} \cdot \frac{1}{r\omega_1(\hat t_{\eta ,H}(r))}.
\end{align}
In order to show (\ref{A37}), we need to slightly adapt the proof of \Cref{T1} by replacing $ir$ with $re^{i\vartheta}$ in (\ref{P1}) and taking into account the constants provided in 
%\todo{correct reference} 
\cite[Theorem 1.1]{langer.pruckner.woracek:heniest}. Then (\ref{A38}) follows as in Step 3 of the proof. \\
For $\vartheta =\frac{\pi}{2}$, the optimal choice of $\eta$ is around $0.13833$ which gives
\[
c_+(0.13833, \frac{\pi}{2}) \approx 1.568, \quad c_-(0.13833, \frac{\pi}{2}) \approx 0.002, \quad \frac{c_+(0.13833, \frac{\pi}{2})}{c_-(0.13833, \frac{\pi}{2})} \approx 675.772 .
\]
While it is possible to derive explicit constants also for $\eta \geq 1-\frac{1}{\sqrt{2}}$, doing so does not result in an improvement of the quotient $c_+ / c_-$.
\end{remark}

\subsection*{Immediate consequences of \Cref{T1}}

\noindent \textit{In order to simplify calculations, unless specified otherwise, we will always assume that $\mr t(r)$ and $\hat t(r)$ are defined implicitly by}
\begin{align}
\label{A12}
(\omega_1 \omega_2)(\mr t(r))=\frac{1}{r^2}, \quad \det \Omega (\hat t(r))=\frac{1}{r^2},
\end{align}
and similarly for $\mr r$ and $\hat r$ (cf. \Cref{A35} with $\eta=2$).
\vspace{1pt}

\noindent We revisit the example from the introduction in more generality. The following example was communicated by Matthias Langer. The calculations can be found in the appendix. \\
\begin{example}
\label{A24}

Let $\alpha > 0$ and $\beta_1, \beta_2 \in \bb R$ where $\beta_1 \neq \beta_2$. Set $\beta_3 := \frac{\beta_1 + \beta_2}{2}$ and define, for $t \in (0,\infty)$,
\begin{align*}
H(t)=
t^{\alpha -1}\left(\begin{matrix}
|\log t|^{\beta_1} & |\log t|^{\beta_3} \\
|\log t|^{\beta_3} & |\log t|^{\beta_2} \\
\end{matrix} \right).
\end{align*}
Then for $r \to \infty$, we have
\begin{itemize}
\item[] $L(r) \asymp (\log r)^{\frac{\beta_1-\beta_2}{2}-2}$ and
\item[] $A(r) \asymp |q_H(ir)| \asymp (\log r)^{\frac{\beta_1-\beta_2}{2}}$,
\end{itemize}
i.e., $L(r) = \Smallo ( A(r))$. Using \Cref{T1}, we can now continue the calculations, leading to
\[
\IM q_H(ir) \asymp (\log r)^{\frac{\beta_1-\beta_2}{2}-1} \asymp \sqrt{L(r)A(r)}.
\]
\end{example}
\vspace{2ex}

\noindent It is an immediate consequence of \Cref{T1} that $\IM q_H$ depends monotonically on the off-diagonal of $H$.
\begin{corollary}
\label{T1+}
Let $H=\smmatrix {h_1}{h_3}{h_3}{h_2}$ and $\tilde{H}=\smmatrix {h_1}{\tilde h_3}{\tilde h_3}{h_2}$ be two Hamiltonians on $[a,b)$. If $t>\hat a (H)$ such that
\[
\Big|\int_a^t h_3(s) \DD s \Big| \geq \Big|\int_a^t \tilde h_3(s) \DD s \Big|,
\]
then
\[
\IM q_H(i \hat r_H(t)) \lesssim \IM q_{\tilde H}(i \hat r_H(t))
\]
with a constant independent of $t$, $H$, and $\tilde H$.
\end{corollary}
\begin{proof}
Our condition states that $|\omega_3(t)| \geq |\tilde \omega_3(t)|$. Taking into account that $t>\hat a(H)$, this means that $0<\det \Omega (t) \leq \det \tilde \Omega (t)$. Hence $\hat r_H(t) \geq \hat r_{\tilde H}(t)$, and further $\hat t_{\tilde H}(\hat r_H(t)) \leq t$. Now, by (\ref{A2}),
\[
\IM q_H(i\hat r_H(t)) \asymp \frac{1}{\hat r_H(t)\omega_2(t)} \leq \frac{1}{\hat r_H(t)\omega_2(\hat t_{\tilde H}(\hat r_H(t)))} \asymp \IM q_{\tilde H}(i\hat r_H(t)).
\]
\end{proof}

\noindent The following result elaborates on the relative behavior of $\IM q_H$ and $|q_H|$. We obtain a quantitative and pointwise relation between $\frac{\IM q_H}{|q_H|}$ and $\frac{\det \Omega}{\omega_1 \omega_2}$, leading to the equivalence
\begin{align}
\label{A14}
\lim_{r \to \infty} \frac{\IM q_H(ir)}{|q_H(ir)|}=0 \,\, \Longleftrightarrow \,\, \lim_{t \to \hat a} \frac{\det \Omega (t)}{(\omega_1 \omega_2)(t)}=0.
\end{align}
The relation between $\frac{\det \Omega}{\omega_1 \omega_2}$ and $\frac{\IM q_H(ir)}{|q_H(ir)|}$ has been investigated also in \cite{langer.pruckner.woracek:gapsatz-arXiv}. Their proof of (\ref{A14})\footnote{In \cite{langer.pruckner.woracek:gapsatz-arXiv}, $\lim_{t \to a}$ was considered instead of $\lim_{t \to \hat a}$.} is based on compactness arguments. \\
Note that our result shows that (\ref{A14}) holds true for $r \to 0$ and $t \to b$ as well.

\begin{proposition}
\label{A4}
Let $H$ be a Hamiltonian on $[a,b)$. Then\footnote{$\mr r(\hat t(r))$ is well-defined because of $\hat t(r) \in (\hat a,b) \subseteq (\mr a,b)$.}
\begin{align}
\label{A9}
\frac{\IM q_H(ir)}{|q_H(ir)|} \asymp \frac{\mr r(\hat t(r))}{r} = \sqrt{\frac{\det \Omega (\hat t(r))}{(\omega_1 \omega_2)(\hat t(r))}} 
\end{align}
for $r \in (0,\infty)$. Moreover,
\begin{align}
\label{A10}
\big|q_H \big( i \mr r(\hat t(r)) \big) \big| \asymp |q_H(ir)|, \quad \quad r \in (0,\infty).
\end{align}
All constants implicit in $\asymp$ do not depend on $H$.
\end{proposition}
\begin{proof}
By definition of $\mr r$ and using (\ref{A2}) and (\ref{A3}),
\[
\mr r(\hat t(r)) = \frac{1}{\sqrt{(\omega_1 \omega_2)(\hat t(r))}} \asymp r \frac{\IM q_H(ir)}{|q_H(ir)|}.
\]
We also have
\[
\sqrt{\frac{\det \Omega (\hat t(r))}{(\omega_1 \omega_2)(\hat t(r))}} = \frac{1}{\sqrt{r^2(\omega_1 \omega_2)(\hat t(r))}} = \frac{\mr r(\hat t(r))}{r},
\]
and (\ref{A9}) follows. \\
For the proof of (\ref{A10}), we need the formula
\[
\omega_1(\mr t(r)) \asymp \frac{|q_H(ir)|}{r}
\]
which we get from \Cref{Y98} applied to $J^{\top}HJ$. Combine this with (\ref{Y35}) to get
\[
|q_H(ir)|^2 \asymp \frac{\omega_1(\mr t(r))}{\omega_2(\mr t(r))}.
\]
On the other hand, (\ref{A2}) and (\ref{A3}) give
\[
|q_H(ir)|^2 \asymp \frac{\omega_1(\hat t(r))}{\omega_2(\hat t(r))}=\frac{\omega_1 \Big(\mr t \big( \mr r(\hat t(r)\big)\Big)}
{\omega_2 \Big(\mr t \big( \mr r(\hat t(r)\big)\Big)} \asymp \big|q_H \big( i \mr r(\hat t(r)) \big) \big|^2.
\]
\end{proof} 

\noindent The freedom in the choice of $\eta$ leads to the following formula that we will refer to later on.

\begin{corollary}
\label{T5}
Let $H$ be a Hamiltonian on $[a,b)$. Then, for any $k>0$,
\begin{align}
\label{A21}
\IM q_H(ikr) \asymp \bigg|q_H(ikr)-\frac{\omega_3(\hat t(r))}{\omega_2(\hat t(r))} \bigg| &\asymp \bigg|q_H(ir)-\frac{\omega_3(\hat t(r))}{\omega_2(\hat t(r))} \bigg|
\end{align}
with constants depending on $k$, but not on $H$. \\
If $\IM q_H(ir)=\Smallo (|q_H(ir)|)$ for $r \to \infty$ \emph{[}$r \to 0$\emph{]}, then
\begin{align}
\label{A20}
q_H(ikr) \sim \frac{\omega_3(\hat t(r))}{\omega_2(\hat t(r))}, \quad \quad r \to \infty \quad [r \to 0].
\end{align}
\end{corollary}
\begin{proof}
Apply \Cref{T1} to $H$ using $\hat t_{1,H}$, and to $kH$ using $\hat t_{k,kH}$. Then $\hat t_{1,H}(r)=\hat t_{k,kH}(r)$, and we write $\hat t(r)$ for short. Keeping in mind that $q_{kH}(z)=q_H(kz)$, this leads to
\[
\IM q_H(ir) \asymp \bigg|q_H(ir)-\frac{\omega_3^{(H)}(\hat t(r))}{\omega_2^{(H)}(\hat t(r))} \bigg| \asymp \frac{1}{r\omega_2^{(H)}(\hat t(r))}
\]
as well as 
\[
\IM q_H(ikr) \asymp \bigg|q_H(ikr)-\frac{k\omega_3^{(H)}(\hat t(r))}{k\omega_2^{(H)}(\hat t(r))} \bigg| \asymp \frac{1}{kr \cdot \omega_2^{(H)}(\hat t(r))}.
\]
(\ref{A21}) follows. Now (\ref{A20}) is obtained by dividing (\ref{A21}) by $|q_H(ikr)|$.
\end{proof}

\section{Behavior of tails of the spectral measure}
\label{S5}

\Cref{T1} that approximately determines the imaginary part of $q_H(ir)$ 
allows us to determine the growth of the spectral measure $\mu_H$ relative 
to suitable comparison functions. Let us introduce the measure $\tilde\mu_H$ on $[0,\infty)$ by
\begin{equation}\label{Y140}
		\tilde\mu_H([0,r)) := \tilde\mu_H(r) := \mu_H((-r,r)), \quad \quad r>0.
\end{equation}
In \Cref{Y190}, equivalent conditions are given for when the function $r \mapsto \tilde \mu_H$ is integrable w.r.t. a given weight function, and also when the measure $\tilde \mu_H$ is finite w.r.t. to a rescaling function. \\
On the other hand, we can view $\tilde\mu_H$ as a function of the positive real parameter $r$, and compare this to a given function $\ms g$. This is what we do in \Cref{Y02}. \\
We note that the content of this section is analogous to \cite[Section 4]{langer.pruckner.woracek:heniest}. The availability of formula (\ref{A2}) leads to improved results in the present article, however we provide less detail as was given in \cite{langer.pruckner.woracek:heniest}.
\newline

\noindent The proofs in this section are based on standard theorems of Abelian-Tauberian type, relating $\mu_H$ to its Poisson integral
\begin{align}
\mc P [\mu_H](z):= \int_{\bb R} \IM \Big( \frac{1}{t-z} \Big) \DD \mu_H(t).
\end{align}
By (\ref{A17}), we have $\mc P [\mu_H](z) = \IM q_H(z) - \beta \IM z $. If $\beta =0$, we can proceed with the application of Abelian-Tauberian theorems without problems. The case $\beta >0$ is equivalent to $a$ being the left endpoint of an $H$-indivisible interval of type $\frac{\pi}{2}$, i.e., $\mr a(H) >a$ and $h_2$ vanishes a.e. on $[a,\mr a(H))$. The restricted Hamiltonian $H_-:=H\big|_{[\mr a(H),b)}$ then has the Weyl coefficient $q_{H_-}(z)=q_H(z)-\beta z$ and thus $\IM q_{H_-}(z) = \mc P [\mu_H](z)$. Hence, we can investigate $\mu_H$ by applying the theorems from this section to $H_-$.

\subsection[{Finiteness of the spectral measure w.r.t. given weight functions}]{Finiteness of the spectral measure w.r.t. given weight functions}
\label{Y190}

\begin{theorem}
\label{AT0}
	Let $H$ be a Hamiltonian defined on $[a,b)$, and assume that $h_2$ does not vanish identically in a neighborhood of $a$.	Let $\ms f$ be a continuous, non-decreasing function,
	and denote by $\mu_H$ the spectral measure of $H$.

	\noindent Then the following statements are equivalent:
	\begin{Enumerate}
	\item
		\begin{equation}
		\label{Y161}
			\int_1^{\infty} \tilde\mu_H(r)\frac{\ms f(r)}{r^3}\DD r<\infty;
		\end{equation}		
	\item
		There is $ a' \in (\hat a,b)$ such that
		\[
			\int_{\hat a}^{a'}
			\frac{1}{\omega_2(t)^2}\binom{\omega_2(t)}{-\omega_3(t)}^*H(t)\binom{\omega_2(t)}{-\omega_3(t)}
			\cdot\ms f\bigl(\det \Omega (t)^{-\frac12}\bigr)\DD t<\infty.
		\]
	\end{Enumerate}

	\noindent If, in addition, $\ms f$ is differentiable,
	then the above conditions hold if and only if there is $ a'\in(\hat a,b)$ such that
	\begin{align*}
		\int_{\hat a}^{a'}
		\frac{(\det \Omega)'(t)}{\omega_2(t)\det \Omega (t)^{\frac 12}} \ms f'\bigl(\det \Omega (t)^{-\frac12}\bigr)\DD t < \infty.
	\end{align*}
\end{theorem}

\begin{proof}
	First note that finiteness of the integrals in the proposition clearly
	does not depend on $a'\in(\hat a,b)$.

\noindent Let $\xi$ be the measure on $[1,\infty)$ such that $\ms f(r)=\xi([1,r))$, $r\ge1$.
	It follows from \cite[Lemma~4]{kac:1982} that
	\[
		\int_{[1,\infty)}\frac{\Poi{\mu_H}(ir)}{r}\DD\xi(r) < \infty
		\quad\Longleftrightarrow\quad
		\int_1^\infty \frac{\tilde\mu_H(r)\ms f(r)}{r^3}\DD r < \infty.
	\]
	Since $h_2$ does not vanish identically in a neighborhood of $a$, we have $\Poi{\mu_H}=\IM q_H$. By \Cref{T1}, we have
	\[
		\frac{\Poi{\mu_H}(ir)}{r}
		\asymp \frac{1}{r^2 \omega_2(\hat t(r))}
		\asymp \frac{\det \Omega (\hat t(r))}{\omega_2(\hat t(r))}.
	\]
Hence
%\todo[color=yellow]{change to eqref}
	\begin{equation}\label{Y162}
		\int_1^{\infty} \tilde\mu_H(r)\frac{\ms f(r)}{r^3}\DD r<\infty
		\;\;\Longleftrightarrow\;\;
		\int_{[1,\infty)}\frac{\det \Omega\bigl(\hat t(r)\bigr)}{\omega_2\bigl(\hat t(r)\bigr)}\DD\xi(r)
		< \infty.
	\end{equation}
	We define a measure $\nu$ on $(0,\infty)$ via $\nu((r,\infty))=\frac{\det \Omega (\hat t(r))}{\omega_2(\hat t(r))}$, $r>0$. Let $\hat\nu$ be the measure on $(\hat a,b)$ satisfying	$\hat\nu((\hat a,t))=\nu((\hat r(t),\infty))=\frac{\det \Omega (t)}{\omega_2(t)}$, $t>\hat a$.
	Integrating by parts (see, e.g., \cite[Lemma~2]{kac:1965}), we can rewrite the first integral in (\ref{Y162}) as follows:
	\begin{align*}
		& \int_{[1,\infty)}\frac{\det \Omega (\hat t(r))}{\omega_2(\hat t(r))} \DD\xi(r)
		= \int_{[1,\infty)}\nu\bigl((r,\infty)\bigr)\DD\xi(r)
		\\
		&= \int_{[1,\infty)}\!\ms f(r)\DD\nu(r) = \int_{(\hat a,\hat t(1)]}\!\ms f(\hat r(t))\DD \hat\nu(t)
		= \int_{(\hat a,\hat t(1)]}\ms f(\hat r(t))\DD \bigg(\frac{\det \Omega}{\omega_2} \bigg)(t)
		\\
		&= \int_{\hat a}^{\hat t(1)}\ms f\bigl(\hat r(t)\bigr)\cdot \frac{1}{\omega_2(t)^2}\binom{\omega_2(t)}{-\omega_3(t)}^*H(t)\binom{\omega_2(t)}{-\omega_3(t)}\DD t.
	\end{align*}
	To prove the additional statement, let us assume that $\ms f$ is differentiable.  Using a substitution we can rewrite
	the second integral in \eqref{Y162} differently:
	\begin{align*}
		&\int_{[1,\infty)}\frac{\det \Omega\bigl(\hat t(r)\bigr)}{\omega_2\bigl(\hat t(r)\bigr)}\DD\xi(r)
		= \int_1^\infty\frac{\det \Omega\bigl(\hat t(r)\bigr)}{\omega_2\bigl(\hat t(r)\bigr)} \ms f'(r)\DD r
		\\
		&= \int_{\hat t(r)}^{\hat a} \!\frac{\det \Omega (t)}{\omega_2(t)}\ms f'(\hat r(t))\hat r'(t)\DD t
		= \frac 12 \int_{\hat a}^{\hat t(r)}\frac{\det \Omega (t)}{\omega_2(t)}\ms f'(\hat r(t))
		\frac{(\det \Omega)'(t)}{\det \Omega (t) ^{\frac32}}\DD t.
	\end{align*}
\end{proof}

\noindent
The following result provides, in particular, information on when the measure $\tilde \mu_H$ is finite w.r.t. a regularly varying rescaling function $\ms g$.

\begin{corollary}
%\label{Y86}
	Let $H$ be a Hamiltonian on $[a,b)$, and assume that $h_2$ does not vanish identically in a neighborhood of $a$.
	Let $\ms g$ be a continuous function that is regularly varying with index $\alpha \in [0,2]$, and denote by $\mu_H$ the spectral measure of $H$ as in (\ref{A17}).
	Then, for $\alpha \in (0,2)$ and every $a'\in (\hat a,b)$, the following statements are equivalent:
\begin{alignat*}{2}
	&\rm (i)\, &&	\int_{[1,\infty)} \frac{\DD\tilde\mu_H(r)}{\ms g(r)} < \infty; \\[1.5ex]
	&\rm (ii)\, &&\int_{\hat a}^{a'}\mkern-10mu
			\frac{1}{\omega_2(t)^2}\binom{\omega_2(t)}{-\omega_3(t)}^* H(t)\binom{\omega_2(t)}{-\omega_3(t)}
			\frac{\DD t}{\det \Omega (t)\ms g\bigl(\det \Omega (t)^{-\frac 12}\bigr)} < \infty . \\[1.5ex]		
	&\rm (iii)\, &&\int_{\hat a}^{a'}\frac{(\det \Omega)'(t)}{\omega_2(t)\det \Omega (t)\ms g\bigl(\det \Omega (t)^{-\frac 12}\bigr)}\DD t
			< \infty;
\end{alignat*}
If $\alpha =0$, then $(iii) \Rightarrow (i)$ and $(iii) \Leftrightarrow (ii)$, while for $\alpha = 2$ we have $(iii) \Rightarrow (i)$ and $(iii) \Rightarrow (ii)$.

%	\begin{alignat*}{5}		
%		&{\r\Omega (i)} \quad&&
%		\begin{array}{c}
%			\scalebox{0.6}{if $\alpha\in(0,2)$} \\[-1.1ex]
%			\scalebox{1.0}{$\Longleftrightarrow$} \\[-0.8ex]
%			\scalebox{1.0}{$\Longrightarrow$} \\[-1.7ex]
%			\scalebox{0.6}{if $\alpha\in\{0,2\}$}
%		\end{array}
%		\quad && {\r\Omega (ia)}
%		\\[-1.5ex]
%		&\Updownarrow
%		\\[-1.5ex]
%		&{\r\Omega (ii)} \quad &&\quad \Longleftrightarrow\quad && {\r\Omega (iia)} \quad&&
%		\begin{array}{c}
%			\scalebox{0.6}{if $\alpha<2$} \\[-1.1ex]
%			\scalebox{1.0}{$\Longleftrightarrow$} \\[-0.8ex]
%			\scalebox{1.0}{$\Longrightarrow$} \\[-1.7ex]
%			\scalebox{0.6}{if $\alpha=2$}
%		\end{array}
%		\quad && {\r\Omega (iib)}.
%	\end{alignat*}
\end{corollary}

\begin{proof}
The increasing function $\ms f(r):=\int_1^r \frac{t}{\ms g(t)} \DD t$ is regularly varying by Karamata's Theorem (\cite[Propositions 1.5.8 and 1.5.9a]{bingham.goldie.teugels:1989}. Moreover,
\begin{equation}
\label{Y137}
	\ms f(r)\;
	\begin{cases}
		\; \asymp\frac{r^2}{\ms g(r)}, & 0 \leq \alpha<2,
		\\[2ex]
		\; \gg\frac{r^2}{\ms g(r)}, & \alpha=2.
	\end{cases}
\end{equation}
Clearly $(iii)$ is equivalent to
\[
\int_{\hat a}^{a'}
		\frac{(\det \Omega)'(t)}{\omega_2(t)\det \Omega (t)^{\frac 12}} \ms f'\bigl(\det \Omega (t)^{-\frac12}\bigr)\DD t < \infty
\]
which is the term appearing in the additional statement of \Cref{AT0}.
Applying \Cref{AT0} and using (\ref{Y137}), this is equivalent to (for $\alpha \in [0,2)$) or implies (for $\alpha =2$) both $(ii)$ and
\[
\int_1^{\infty} \tilde\mu_H(r)\frac{dr}{r\ms g(r)}<\infty.
\]
By \cite[Proposition 4.5]{langer.pruckner.woracek:heniest}, this is further equivalent to (for $\alpha \in (0,2]$) or implies (for $\alpha =0$) the first item.
\end{proof}

\subsection{Comparative growth of the distribution function}
\label{Y02}

In this section we investigate $\limsup$-conditions for the
quotient $\frac{\tilde\mu_H(r)}{\ms g(r)}$ instead of integrability conditions.
Let us introduce the corresponding classes of measures.

\begin{definition}\label{Y77}
	Let $\ms g(r)$ be a regularly varying function with index $\alpha \in [0,2]$
	and $\lim_{r\to\infty}\ms g(r)=\infty$.
	Then we set
	\[
		\mc F_{\ms g} \DE \big\{\mu\DS \tilde\mu(r) \lesssim \ms g(r), \, r \to \infty \big\},\qquad
		\mc F_{\ms g}^0 \DE \big\{\mu\DS \tilde\mu(r) = \Smallo ( \ms g(r)), \, r \to \infty\big\},
	\]
	where again $\tilde\mu(r)\DE\mu((-r,r))$.
\end{definition}

\noindent It should be mentioned that, for non-decreasing $\ms g$, if 
\[
	\int_{[1,\infty)} \frac{\DD\tilde\mu(r)}{\ms g(r)} < \infty ,
\]
then $\mu \in \mc F_{\ms g}^0 \subseteq \mc F_{\ms g}$. For further discussion of this relation, the reader is referred to \cite{langer.pruckner.woracek:heniest}.

\begin{theorem}\label{Y74}
	Let $H$ be a Hamiltonian on $[a,b)$, and assume that $h_2$ does not vanish identically in a neighborhood of $a$.
	Let $\ms g(r)$ be a regularly varying function with index $\alpha \in [0,2]$ and $\lim_{r\to\infty}\ms g(r)=\infty$. Denote by $\mu_H$ the spectral measure of $H$. For $\alpha < 2$, the following statements hold:

\begin{alignat*}{4}
	&\rm (i)\quad &&\mu_H \in \mathcal{F}_{\ms g} \quad && 
	\Leftrightarrow
		&& \quad \limsup_{t \to \hat a} \frac{1}{\omega_2(t)\ms g \big( \det \Omega (t)^{-\frac 12} \big)}<\infty; \\[1ex]		
		&\rm (ii)\quad &&\mu_H \in \mathcal{F}_{\ms g}^0 \quad &&
	\Leftrightarrow
		&&\quad \lim_{t \to \hat a} \frac{1}{\omega_2(t)\ms g \big( \det \Omega (t)^{-\frac 12} \big)}=0.	
\end{alignat*}
If $\alpha =2$, then the right hand side of $(i)$, $(ii)$ implies the left hand side, respectively.
\end{theorem}

\begin{proof}
	We use \cite[Lemma 4.16]{langer.pruckner.woracek:heniest} which, adapted to our situation, reads as
	\[
		c_{\alpha}\limsup_{r\to\infty}\biggl(\frac{r}{\ms g(r)}\Poi{\mu_H}(ir)\biggr) \leq \limsup_{r\to\infty}\frac{\tilde\mu_H(r)}{\ms g(r)}
		\leq c_{\alpha}' \limsup_{r\to\infty}\biggl(\frac{r}{\ms g(r)}\Poi{\mu_H}(ir)\biggr),
	\]
	and the second inequality holds even for $\alpha =2$. Since $h_2$ does not vanish identically in a neighborhood of $a$, we have $\Poi{\mu_H}=\IM q_H$. Therefore, the assertion follows from \Cref{T1} and a substitution $r=\hat r(t)$.
\end{proof}

\section{Weyl coefficients with tangential behavior}
\label{S3}

In this section, we investigate the scenario 
\begin{align}
\label{A16}
\lim_{r \to \infty} \frac{\IM q_H(ir)}{|q_H(ir)|}=0 \quad \quad \text{or} \quad \quad \liminf_{r \to \infty} \frac{\IM q_H(ir)}{|q_H(ir)|}=0.
\end{align}
This is equivalent to tangential behavior of $q_H(ir)$, i.e.,
\[
\lim_{r \to \infty} \arg q_H(ir) \in \{0,\pi\} \quad \text{or} \quad \liminf_{r \to \infty} \min \big\{\arg q_H(ir),\, \pi-\arg q_H(ir) \big\}=0.
\]

%\begin{remark} 
%\label{R1}
%In the situation of (\ref{A16}), our \Cref{T1} is always applicable. To see this, recall that $(a,b)$ starts with an $H$-indivisible interval if and only if $\IM q_H(ir) \asymp \frac 1r$ or $\IM q_H(ir) \asymp r$ for $r \to \infty$, see e.g. \cite[Theorems 4.33, 4.34]{remling:2018}. In this situation, we must have $\liminf_{r \to \infty} \frac{\IM q_H(ir)}{|q_H(ir)|} > 0$, which follows from the integral representation (\ref{A17}).
%\end{remark}

%\vspace{1.5ex}

From \Cref{A4} we get that
\begin{align}
\label{A25}
\lim_{n \to \infty} \frac{\IM q_H(i r_n)}{|q_H(i r_n)|}=0 \,\, \Longleftrightarrow \,\, \lim_{n \to \infty} \frac{\det \Omega(\hat t(r_n))}{(\omega_1 \omega_2)(\hat t(r_n))}=0.
\end{align}
for every sequence $r_n \to \infty$. 
All results in this section can be seen from the canonical systems perspective as well as from the Herglotz functions perspective. \\

\vspace{1pt}

\noindent To start with, we observe that the second assertion in (\ref{A16}) implies the first unless the limit inferior is assumed only along very sparse sequences. We formulate this fact in the language of Herglotz functions, and prove it within the canonical systems setting. However, we do not know a purely function theoretic proof (which may very well exist in the literature).

\begin{lemma}
Let $q$ be a Herglotz function. Suppose there is a sequence $(r_n)_{n \in \bb N}$ with $r_n \to \infty$, $\sup_{n \in \bb N} \frac{r_{n+1}}{r_n} < \infty$, and 
\[
\lim_{n \to \infty} \frac{\IM q(ir_n)}{|q(ir_n)|}=0.
\]
Then $\lim_{r \to \infty} \frac{\IM q(ir)}{|q(ir)|}=0$.
\end{lemma}
\begin{proof}
Let $H$ be a Hamiltonian (on $[0,\infty)$), such that $q=q_H$. \\
Let $d(t):=\frac{\det \Omega(t)}{(\omega_1 \omega_2)(t)}$. Set $t_n := \hat t(r_n)$, then by (\ref{A9}),
\[
d(t_n) \asymp \bigg( \frac{\IM q(ir_n)}{|q(ir_n)|} \bigg)^2 \xrightarrow{n \to \infty} 0.
\]
Suppose that the assertion was not true, i.e., there is a sequence $\xi_1 > \xi_2 > ...$ converging to $0$, such that $d(\xi_k) \geq C>0$ for all $k$. For $k \in \bb N$, set $n(k):=\max \{n \in \bb N \mid t_n > \xi_k \}$. We obtain
\begin{align*}
&\Big(\frac{r_{n(k)+1}}{r_{n(k)}} \Big)^2 = \frac{\det \Omega(t_{n(k)})}{\det \Omega (t_{n(k)+1})} \geq \frac{\det \Omega(\xi_k)}{\det \Omega (t_{n(k)+1})} \\
&=\frac{d(\xi_k)}{d(t_{n(k)+1})} \cdot \frac{(\omega_1 \omega_2)(\xi_k)}{(\omega_1 \omega_2)(t_{n(k)+1})} \geq \frac{C}{d(t_{n(k)+1})} \xrightarrow{k \to \infty} \infty
\end{align*}
which contradicts our assumption.
\end{proof}

\noindent Recall formulae (\ref{A10}) and (\ref{A9}). On an intuitive level, they tell us that in the case that $\IM q_H(ir) \not\asymp |q_H(ir)|$, the growth of $|q_H(ir)|$ is restricted since $\mr r(\hat t(r))$ is then far away from $r$. If read in the other direction, this means that if $|q_H(ir)|$ grows quickly and without oscillating too much, then $\mr r(\hat t(r))$ and $r$ should be close to each other, and hence the quotient $\frac{\IM q_H(ir)}{|q_H(ir)|}$ should not decay.\\
The following definition introduces the notions needed in Theorems \ref{T7} and \ref{T8}, which confirm this intuition.

\begin{definition}
\label{A48}
\item[$\rhd$] A measurable function $f: (0,\infty) \to (0,\infty)$ is called \textit{regularly varying (at infinity) with index $\alpha \in \bb R$} if, for any $\lambda >0$,
\begin{align}
\lim_{r \to \infty}\frac{f(\lambda r)}{f(r)} = \lambda^{\alpha}.
\end{align}
If $\alpha =0$, then $f$ is also called \textit{slowly varying (at infinity)}.
\item[$\rhd$] A measurable function $f: (0,\infty) \to (0,\infty)$ is \textit{positively increasing (at infinity)} if there is $\lambda \in (0,1)$ such that
\begin{align}
\limsup_{r \to \infty} \frac{f(\lambda r)}{f(r)} <1.
\end{align}
Let us say explicitly that we do not require $f$ to be monotone.
\end{definition}

\stepcounter{lemma}
\begin{subtheorem}
\label{T7}
Let $q \neq 0$ be a Herglotz function. If $|q(ir)|$ or $\frac{1}{|q(ir)|}$ is positively increasing at infinity (in particular, if $|q(ir)|$ is regularly varying with index $\alpha \neq 0$), then $\IM q(ir) \asymp |q(ir)|$ as $r \to \infty$.
\end{subtheorem}

\begin{subtheorem}
\label{T8}
Let $q \neq 0$ be a Herglotz function. If $\IM q(ir) = \Smallo (|q(ir)|)$ as $r \to \infty$, then, for every $\delta \in [0,1)$,
\begin{align}
\label{A19}
\lim_{r \to \infty} \frac{ q \Big(ir \Big[\frac{\IM q(ir)}{|q(ir)|}\Big]^{\delta} \Big)}{q(ir)} =1.
\end{align}
For $k>0$, we also have $\lim_{r \to \infty} \frac{q(ikr)}{q(ir)}=1$, in particular, $|q(ir)|$ is slowly varying at infinity.
\end{subtheorem}

\begin{remark}
%\Cref{T7} raises the following questions: Is the assumption that $|q(ir)|$ or $\frac{1}{|q(ir)|}$ be positively increasing vital? Or could it be replaced by, say, $|q(ir)| \gtrsim r^{\delta}$ for $r \to \infty$, for some $\delta > 0$? \\
In \Cref{T7}, the requirement that $|q(ir)|$ should be positively increasing is meaningful. It is not enough that $|q(ir)|$ grows sufficiently fast, say, $|q(ir)| \gtrsim r^{\delta}$ for $r \to \infty$ and some $\delta > 0$. \\
In fact, for any given $\delta \in (0,1)$, we construct in \Cref{T6} a Hamiltonian\footnote{Choose suitable parameters $p,l \in (0,1)$, such that $\delta = \frac{\log l}{\log (pl)}$, i.e., $p=l^{\delta ^{-1}-1}$.} $H$ whose Weyl coefficient $q_H$ satisfies (see \Cref{R2}) $|q_H(ir)| \gtrsim r^{\delta}$ as $r \to \infty$, but
\[
\liminf_{r \to \infty} \frac{\IM q_H(ir)}{|q_H(ir)|} = 0.
\]
In other words, $\IM q_H(ir) \not\asymp |q_H(ir)|$. \\
Note also that for the above-mentioned $H$, certainly $|q_H(ir)|$ is not slowly varying \cite[Proposition 1.3.6]{bingham.goldie.teugels:1989}. Hence, in \Cref{T8} it is not enough to require $\IM q(ir_n) = \Smallo (|q(ir_n)|)$ on some sequence $r_n \to \infty$.
\end{remark}

\begin{example}
Let $q(z)=\log z$, satisfying $|q(ir)| = \big[(\log r)^2+\frac{\pi^2}{4}\big]^{1/2}$ which is increasing. However, $\IM q(ir)$ is constant and hence $\IM q(ir) = \Smallo (|q(ir)|)$ as $r \to \infty$. \Cref{T7} fails because $|q(ir)|$ is not positively increasing.
\end{example}

\begin{proof}[Proof of \Cref{T7}]
Assume first that $|q(ir)|$ is positively increasing. Then there are $\lambda, \sigma \in (0,1)$ and $R > 0$ such that
\begin{align}
\label{A18}
\frac{|q(i\lambda r)|}{|q(ir)|} \leq \sigma, \quad \quad r \geq R.
\end{align}
Let $H$ be a Hamiltonian with Weyl coefficient $q_H=q$, allowing us to use (\ref{A10}). \\
Suppose that the assertion was not true. Then there is a (w.l.o.g., monotone) sequence $r_n \to \infty$ with $\lim_{n \to \infty} \frac{\IM q(ir_n)}{|q(ir_n)|}=0$. Let $m(n)$ be such that 
\[
\lambda^{m(n)+1} \leq \frac{\mr r(\hat t(r_n))}{r_n} < \lambda^{m(n)}.
\]
Note that $m(n) \to \infty$ because of (\ref{A9}). \\
Furthermore, (\ref{A10}) ensures that there is $\beta >0$ with
\[
\beta \leq \frac{|q(i \mr r(\hat t(r)))|}{|q(ir)|}, \quad r \in (0,\infty).
\]
We will also need that for $0<r<r'$,
\[
\frac{|q(ir)|}{|q(ir')|} \asymp \frac{r' \omega_2(\mr t(r'))}{r \omega_2(\mr t(r))} \leq \frac{r'}{r}
\]
because $\omega_2$ is nondecreasing. \\
Choosing $n$ so big that $\mr r(\hat t(r_n)) \geq R$, we get the contradiction
\begin{align*}
\beta &\leq \frac{\big|q \big(i \mr r(\hat t(r_n)) \big) \big|}{\big|q \big(ir_n \big)\big|} = \frac{\big|q \big(i \mr r(\hat t(r_n)) \big) \big|}{\big|q \big(i\lambda^{m(n)}r_n \big)\big|} \cdot \prod_{j=0}^{m(n)-1}\frac{\big|q \big(i\lambda^{j+1} r_n \big)\big|}{\big|q \big(i\lambda^j r_n \big)\big|} \\
&\lesssim \frac{\lambda^{m(n)}r_n}{\mr r(\hat t(r_n))} \sigma^{m(n)}
\leq \frac{\sigma^{m(n)}}{\lambda} \xrightarrow{n \to \infty} 0.
\end{align*}
This proves the theorem in the case that $|q(ir)|$ is positively increasing. \\
If, on the other hand, $\frac{1}{|q(ir)|}$ is positively increasing, we may set $\tilde q :=-\frac 1q$, for which $|\tilde q(ir)|$ is positively increasing. We obtain
\[
\frac{\IM q(ir)}{|q(ir)|} = \frac{\IM \tilde q(ir)}{|\tilde q(ir)|} \asymp 1.
\]
Finally, we note that if $|q(ir)|$ is regularly varying with index $\alpha >0$, then it is also positively increasing. If $|q(ir)|$ is regularly varying with index $\alpha < 0$, then $\frac{1}{|q(ir)|}$ is regularly varying with index $-\alpha > 0$ and thus positively increasing.
\end{proof}

\noindent Our proof of \Cref{T8} is elementary - only folklore facts that follow from the Herglotz integral representation (\ref{A17}) are needed. We would be interested in an elementary proof of \Cref{T7} as well, which so far we have not found. \\

\noindent One fact needed in the following proof is the following: For any Herglotz function $q$ and any $z \in \bb C_+$, we have
\begin{equation}
\label{A5}
|q'(z)| \leq \frac{\IM q(z)}{\IM z}.
\end{equation}
This can be seen using the representation (\ref{A17}): We write
\[
q'(z)=b+\int_{\bb R} \frac{d\sigma (t)}{(t-z)^2}
\]
and obtain
\[
|q'(z)| \leq b+\int_{\bb R} \frac{d\sigma (t)}{|t-z|^2}= \frac{\IM q(z)}{\IM z}.
\]

\begin{proof}[Proof of \Cref{T8}]
Let $k \in (0,1)$. Then
\begin{align}
\label{A15}
&|\log q \big(ikr \big)-\log q(ir)|=\Big|\int_{kr}^r i(\log q)'(is) \DD s \Big| \leq \int_{kr}^r \big|(\log q)'(is)\big| \DD s. 
\end{align}
Apply (\ref{A5}) to $\log q$ and to $i\pi-\log q$ to obtain
\begin{align*}
&|(\log q)'(is)| \leq \frac 1s \min \big\{\IM [\log q(is)],\pi - \IM [\log q(is)] \big\} \\
&= \frac 1s \min \big\{\arg q(is), \pi-\arg q(is) \big\} \asymp \frac{\IM q(is)}{s|q(is)|}
\end{align*}
for all $s>0$. We will also need monotonicity in $s$ of $s\frac{\IM q(is)}{|q(is)|}$. In fact, it is easy to see from (\ref{A17}) that $s \IM q(is)$ is nondecreasing in $s$. Now we can write
\[
s\frac{\IM q(is)}{|q(is)|} = \sqrt{s \IM q(is) \cdot s \IM \Big(-\frac{1}{q(is)} \Big)}
\]
and hence $s\frac{\IM q(is)}{|q(is)|}$ is nondecreasing in $s$.
Putting together and continuing the estimation in (\ref{A15}), we obtain
\begin{align}
&|\log q \big(ikr \big)-\log q(ir)| \lesssim \int_{kr}^r \frac{\IM q(is)}{s|q(is)|} \DD s \leq r \frac{\IM q(ir)}{|q(ir)|} \cdot \int_{kr}^r \frac{ds}{s^2} \nonumber \\
&= r \frac{\IM q(ir)}{|q(ir)|}\Big(\frac{1}{kr}-\frac 1r \Big) \asymp \frac{\IM q(ir)}{|q(ir)|} \xrightarrow{r \to \infty} 0. \label{A28}
\end{align}
This shows $\lim_{r \to \infty} \frac{q(ikr)}{q(ir)}=1$. To prove (\ref{A19}), set $k(r):= \frac{\IM q(ir)}{|q(ir)|}$ and repeat the calculations up to the second to last term in (\ref{A28}), but with $k$ replaced by $k(r)^{\delta}$, where $\delta \in [0,1)$. Since
\[
r k(r) \Big(\frac{1}{rk(r)^{\delta}}-\frac 1r \Big) \asymp k(r)^{1-\delta} \xrightarrow{r \to \infty} 0,
\]
we arrive at (\ref{A19}).
\end{proof}
\noindent Note that $\lim_{r \to \infty} \frac{q(ikr)}{q(ir)}=1$ is also a consequence of (\ref{A20}). The preceding proof, in addition to being elementary, is needed to show (\ref{A19}) which, upon taking absolute values, can be seen as slow variation with a rate.

\section{Maximal oscillation within Weyl disks}
\label{S4}

\noindent In order to explain the aim of this section, let us first recall the notion of Weyl disks. Let $W(t,z) \in \bb C^{2 \times 2}$ be the fundamental solution of
\begin{equation}
\frac{d}{dt} W(t,z)J=zW(t,z)H(t),
\end{equation}
with initial condition $W(a,z)=I$, solving the transpose of equation (\ref{A33}). We define the \textit{Weyl disks}
\begin{equation}
\label{A34}
D_{t,z}:=\Big\{ \frac{w_{11}(t,z)\tau + w_{12}(t,z)}{w_{21}(t,z)\tau + w_{22}(t,z)} \Big| \tau \in \overline{\bb C_+} \Big\} \subseteq \overline{\bb C_+},
\end{equation}
where $\bb C_+=\{z \in \bb C | \IM z>0\}$, and the closure is taken in the Riemann sphere $\overline{\bb C}=\bb C \cup \{ \infty \}$. For fixed $z \in \bb C_+$ and $t_1 \leq t_2$, we have $D_{t_1,z} \supseteq D_{t_2,z}$, and the disks shrink down to a single point which is $q_H(z)$:
\[
\bigcap_{t \in [a,b)} D_{t,z} = \{ q_H(z) \}.
\]

\noindent Now we review the estimate (\ref{A43}) which has a geometric interpretation. Namely, the functions $L(r)$ and $A(r)$ give, up to $\asymp$, the imaginary part of the bottom and top point of $D_{\mr t(r),ir}$, respectively. The size of $\IM q_H(ir)$ relative to $L(r)$ and $A(r)$ thus corresponds to the vertical position of $q_H(ir)$ within the disk $D_{\mr t(r),ir}$. \\
In this section we give answers to several questions from \cite{langer.pruckner.woracek:heniest}. For instance, the question was raised whether there is a Hamiltonian $H$ for which $L(r) \asymp \IM q_H(ir) \not\asymp A(r)$ for $r \to \infty$. The answer to this particular question is no, cf. \Cref{T10}. However\footnote{In this section, we use the more transparent notation $f(r) \ll g(r)$ instead of $f(r) = \Smallo (g(r))$.}, $L(r_n) \asymp \IM q_H(ir_n) \ll A(r_n)$ on a subsequence $r_n \to \infty$ is possible, and we provide examples for this in \Cref{T6} and in \Cref{R3}. The Weyl coefficient of the Hamiltonian constructed in \Cref{T6} exhibits "maximal" oscillatory behaviour in the sense that it goes back and forth between the bottoms and tops of the disks $D_{\mr t(r),ir}$.

\begin{proposition}
\label{T10}
Let $H$ be a Hamiltonian on $(a,b)$. The following statements hold:
\begin{itemize}
\item[$(i)$] Suppose that $L(r) \not\asymp A(r)$ as $r \to \infty$. Then there exists a sequence $(r_n)_{n \in \bb N}$ such that $r_n \to \infty$, $L(r_n) \ll A(r_n)$, and
\[
\IM q_H(ir_n) \gtrsim \sqrt{L(r_n)A(r_n)}.
\]
\item[$(ii)$] Suppose that $L(r) \not\asymp A(r)$, but not $L(r) \ll A(r)$ as $r \to \infty$. Then there is also $(r_n')_{n \in \bb N}$ with $r_n' \to \infty$, $L(r_n') \ll A(r_n')$, and
\begin{equation}
\label{A32}
\IM q_H(ir_n') \asymp \sqrt{L(r_n')A(r_n')}.
\end{equation}
\end{itemize}

\end{proposition}
\begin{proof}
We shorten notation by setting $d(t):=\frac{\det \Omega(t)}{(\omega_1 \omega_2)(t)}$. By assumption, $\liminf_{t \to \hat a} d(t)=0$. Let $c \in (\hat a,b)$ be fixed and set $t_n := \max\{t \leq c \mid d(t) \leq \frac 1n \}$. With $t_n^+:= \hat t(\mr r(t_n)) \geq t_n$, we have $d(t_n^+) \geq \frac 1n = d(t_n)$ if $n$ is large enough for $t_n^+ \leq c$ to hold. Using (\ref{A9}), we obtain
\[
\bigg(\frac{\IM q_H(i \mr r(t_n))}{A(\mr r(t_n))} \bigg)^2 \asymp d(t_n^+) \geq d(t_n)=\frac{L(\mr r(t_n))}{A(\mr r(t_n))}.
\]
Note that $L(\mr r(t_n)) \ll A(\mr r(t_n))$ because of $d(t_n) \to 0$. \\
Suppose now that $s:=\limsup_{r \to \infty} \frac{L(r)}{A(r)}>0$. Set $\xi_n := \max \{t \leq t_n \mid d(\xi_n)=\frac s2\}$ and find $\tau_n$ between $\xi_n$ and $t_n$ such that $d(\tau_n) = \min \{d(t) \mid t \in [\xi_n, t_n]\}$. Certainly, $d(\tau_n) \leq d(t_n)=\frac 1n$ and $d(\tau_n) \leq d(t)$ for all $t \in [\xi_n,c]$. Also note that by the same arguments as above,
\begin{align}
\label{A26}
\IM q_H(i \mr r(\tau_n)) \gtrsim \sqrt{L(\mr r(\tau_n)A(\mr r(\tau_n))}.
\end{align}

\noindent We prove next that $\hat r(\tau_n) \ll \mr r(\xi_n)$. Note that by passing to a subsequence and possibly switching signs of $\omega_3$ by looking at $J^{\top}HJ$ instead of $H$, we can assume that 
\[
\lim_{n \to \infty} \frac{\omega_3(\tau_n)}{\sqrt{(\omega_1 \omega_2)(\tau_n)}}=1.
\]
A calculation shows that $\sqrt{(\omega_1 \omega_2)(t)}-\omega_3(t)$ is increasing. Hence
\begin{align}
\label{A27}
&\bigg(\frac{\hat r(\tau_n)}{\mr r(\xi_n)} \bigg)^2= \frac{(\omega_1 \omega_2)(\xi_n)}{\det \Omega(\tau_n)} \\
&=\frac{(\omega_1 \omega_2)(\xi_n)\Big(1-\frac{\omega_3(\tau_n)}{\sqrt{(\omega_1 \omega_2)(\tau_n)}}\Big)}{\Big(\sqrt{(\omega_1 \omega_2)(\tau_n)}-\omega_3(\tau_n)\Big)^2 \Big(1+\frac{\omega_3(\tau_n)}{\sqrt{(\omega_1 \omega_2)(\tau_n)}}\Big)} \nonumber\\
&\leq \frac{(\omega_1 \omega_2)(\xi_n)\Big(1-\frac{\omega_3(\tau_n)}{\sqrt{(\omega_1 \omega_2)(\tau_n)}}\Big)}{\Big(\sqrt{(\omega_1 \omega_2)(\xi_n)}-\omega_3(\xi_n)\Big)^2 \Big(1+\frac{\omega_3(\tau_n)}{\sqrt{(\omega_1 \omega_2)(\tau_n)}}\Big)} \nonumber\\
&= \frac{\Big(1-\frac{\omega_3(\tau_n)}{\sqrt{(\omega_1 \omega_2)(\tau_n)}}\Big)}{\Big(1-\frac{\omega_3(\xi_n)}{\sqrt{(\omega_1 \omega_2)(\xi_n)}}\Big)^2 \Big(1+\frac{\omega_3(\tau_n)}{\sqrt{(\omega_1 \omega_2)(\tau_n)}}\Big)} \lesssim 1-\frac{\omega_3(\tau_n)}{\sqrt{(\omega_1 \omega_2)(\tau_n)}} \to 0. \nonumber
\end{align}
Let $\tau_n^- := \mr t(\hat r(\tau_n))$. By the calculation above, $\mr r(\tau_n^-)=\hat r(\tau_n) < \mr r(\xi_n)$ for large enough $n$, implying $\tau_n^- > \xi_n$ and hence $d(\tau_n^-) \geq d(\tau_n)$. Consequently,
\[
\frac{L(\hat r(\tau_n))}{A(\hat r(\tau_n))} =d(\tau_n^-) \geq d(\tau_n) \asymp \bigg(\frac{\IM q_H(i\hat r(\tau_n))}{A(\hat r(\tau_n))} \bigg)^2.
\]
This means that
\[
\IM q_H(i \hat r(\tau_n)) \leq C \sqrt{L(\hat r(\tau_n)) A(\hat r(\tau_n))}
\]
for some $C>0$ and all large $n$. Recall (\ref{A26}) and choose $C'>0$, w.l.o.g. $C'<C$, such that
\[
\IM q_H(i \mr r(\tau_n)) \geq C' \sqrt{L(\mr r(\tau_n)) A(\mr r(\tau_n))}
\]
for large $n$. By continuity, we find, for each large $n$, an $r_n' \in [\mr r(\tau_n),\hat r(\tau_n)]$ with 
\[
\frac{\IM q_H(ir_n')}{\sqrt{L(r_n')A(r_n')}} \in [C',C],
\]
such that $(r_n')_{n \in \bb N}$ satisfies (\ref{A32}). The only thing left to prove is that $L(r_n') \ll A(r_n')$. \\
Suppose not, then on a subsequence we would have $L(r_n') \asymp A(r_n')$. Consider $\xi_n':=\mr t(r_n') \leq \tau_n$ which would then satisfy $d(\xi_n') \gtrsim 1$ and hence
\[
1-\frac{\omega_3(\xi_n')}{\sqrt{(\omega_1 \omega_2)(\xi_n')}} \gtrsim 1.
\]
Now look at (\ref{A27}), but with $\xi_n$ replaced by $\xi_n'$. It follows that, for large $n$, $\hat r(\tau_n) < r_n'$, contradicting the choice of $r_n'$.
\end{proof}

%\todo{delete this paragraph?}
%\noindent We return to the situation of \Cref{T8}. Certainly, if $|q_H(ir)|$ is not slowly varying (for example, if $|q_H(ir)| \gtrsim r^{\delta}$ for some $\delta>0$), then
%\[
%\limsup_{t \to a} \frac{\det \Omega(t)}{(\omega_1 \omega_2)(t)} >0.
%\]
%In fact, this inference is sharp - we will give an example where the limit inferior is $0$. However, by \Cref{T7}, the example has to be such that $|q_H(ir)|$ is not regularly varying. Our strategy is to prescribe $\frac{\det \Omega(t)}{(\omega_1 \omega_2)(t)}$, forcing oscillatory behaviour.

\noindent In the following definition, we construct a Hamiltonian by prescribing $f:=\frac{\omega_3}{\sqrt{\omega_1 \omega_2}}$ and choosing $f$ to be a highly oscillating function. It should be mentioned that the method we use for prescription works on a general basis: Any locally absolutely continuous function with values in $(-1,1)$ occurs as $\frac{\omega_3}{\sqrt{\omega_1 \omega_2}}$ for some Hamiltonian. Details can be found in the appendix.

\begin{definition}
\label{T6}
Let $(t_n)_{n \in \bb N}, (\xi_n)_{n \in \bb N}$ be sequences of positive numbers converging to zero, where $\xi_{n+1}<t_n<\xi_n$ for all $n \in \bb N$. Choose $p,l \in (0,1)$ and set
\[
f(t_n)=1-p^n, \quad \quad f(\xi_n)=l^n
\]
and interpolate between those points using monotone and absolutely continuous functions (e.g., linear interpolation). Set
\[
\alpha_1(t):= \left\{ \begin{matrix}
\frac{f'(t)}{1-f(t)}, & t \in (\xi_{n+1},t_n), \\
0, & t \in (t_n,\xi_n)
\end{matrix}\right.
\]
and
\[
\alpha_2(t):= \left\{ \begin{matrix}
\frac{f'(t)}{1-f(t)}, & t \in (\xi_{n+1},t_n), \\
-2\frac{f'(t)}{f(t)}, & t \in (t_n,\xi_n)
\end{matrix}\right.
\]
For $t \in [0,t_1]$, let $\omega_i(t):=\exp \Big(-\int_t^{t_1} \alpha_i(s) \DD s \Big)$, $i=1,2$, and $\omega_3(t):=\sqrt{(\omega_1 \omega_2)(t)} \cdot f(t)$. Set $h_i(t)=\omega_i'(t)$, $i=1,2,3$, $t \in [0,t_1]$. For $t \in (t_1,\infty)$, let $h_1(t):=1$ and $h_2(t):=h_3(t):=0$. Finally, define
\[
H_{p,l} :=\begin{pmatrix}
h_1 & h_3 \\
h_3 & h_2
\end{pmatrix}.
\]
\end{definition}

\begin{lemma}
$H_{p,l}$ is a Hamiltonian on $[0,\infty)$, and $\omega_i(t)=\int_0^t h_i(s) \DD s$ for $i=1,2,3$ and $t \in [0,t_1]$. Moreover, $0$ is not the left endpoint of an $H_{p,l}$-indivisible interval.
\end{lemma}
\begin{proof}
We write $H$ instead of $H_{p,l}$ for short. First we show that $H(t) \geq 0$ for all $t \in [0,t_1]$. Start by noting that, for $i=1,2$,
\[
\frac{h_i(t)}{\omega_i(t)}=(\log \omega_i)'(t)=\alpha_i(t),
\]
and calculate
\begin{align*}
&\frac{h_3(t)^2}{(\omega_1 \omega_2)(t)}=\frac{\big[(\sqrt{\omega_1 \omega_2}f )'(t)\big]^2}{(\omega_1 \omega_2)(t)} = \Big(f'(t)+\frac 12 \Big[\frac{h_1(t)}{\omega_1(t)}+\frac{h_2(t)}{\omega_2(t)} \Big]f(t) \Big)^2 \\
&=\Big(f'(t)+\frac{\alpha_1(t)+\alpha_2(t)}{2} f(t) \Big)^2.
\end{align*}
If $t \in (t_n,\xi_n)$, then this equates to $0$, as does
\[
\frac{(h_1h_2)(t)}{(\omega_1 \omega_2)(t)}=\alpha_1(t)\alpha_2(t)=0.
\]
For $t \in (\xi_{n+1},t_n)$,
\[
\Big(f'(t)+\frac{\alpha_1(t)+\alpha_2(t)}{2} f(t) \Big)^2=\Big(\frac{f'(t)}{1-f(t)}\Big)^2=\alpha_1(t)\alpha_2(t)=\frac{(h_1h_2)(t)}{(\omega_1 \omega_2)(t)}.
\]
In both cases, $\det H(t)=0$. For $i=1,2$, as $\alpha_i(t) \geq 0$, $t \in [0,t_1]$, certainly $\omega_i(t)$ is increasing and thus $h_i(t) \geq 0$. This suffices to show that $H(t) \geq 0$. \\
$H$ is in limit point case since, for $t>t_1$, the trace of $H(t)$ equals $1$. To show that $\omega_i(t)=\int_0^t h_i(s) \DD s$, $i=1,2,3$, $t \in [0,t_1]$, we need to check that $\lim_{t \to 0} \omega_i(t)=0$. For $i=1$, this follows from
\begin{equation}
\int_0^{t_1} \alpha_1(s) \DD s = \sum_{n=1}^{\infty} \int_{\xi_{n+1}}^{t_n} \frac{f'(s)}{1-f(s)} \DD s=\sum_{n=1}^{\infty} \big[\log (1-l^{n+1})-\log(p^n) \big]=\infty.
\end{equation}
For $i=2$, it follows from the fact that $\alpha_2(t) \geq \alpha_1(t)$ for all $t \in [0,t_1]$, and for $i=3$ it follows from the definition of $\omega_3$ and the fact that $f(t) <1$, $t \in [0,t_1]$. \\
Finally, $0$ is not the left endpoint of an $H$-indivisible interval because 
\[
\det \Omega(t)=(\omega_1 \omega_2)(t) \big(1-f(t)^2 \big) >0
\]
for all $t \in (0,t_1]$.
\end{proof}

%\noindent From the fact that $\frac{\det \Omega(t)}{(\omega_1 \omega_2)(t)}=1-f(t)^2$ we immediately obtain
%\[
%0=\liminf_{t \to 0} \frac{\det \Omega(t)}{(\omega_1 \omega_2)(t)} < \limsup_{t \to 0} \frac{\det \Omega(t)}{(\omega_1 \omega_2)(t)}=1.
%\]

\noindent We investigate the behaviour for $r \to \infty$ of $\IM q_{H_{p,l}}(ir)$ as well as $L(r)$ and $A(r)$. A rough description of the situation is: 

\begin{center}
\label{tikz1}
\begin{tikzpicture}
\draw[scale=0.5, loosely dashed, domain=1:3.65, smooth, variable=\x] plot ({(\x+2)*(\x+2)}, {4*(\x+2)});

%A(r)
\draw[out=30, in=240] (4.5,6.3) to (6,8.5);
\draw[out=60, in=190] (6,8.5) to (6.7,9);
\draw[out=10, in=175] (6.7,9) to (7.5,8.9);
\draw[out=-5, in=190] (7.5,8.9) to (8,9);
\draw[out=10, in=205] (8,9) to (9.5,9);

\draw[out=25, in=240] (9.5,9) to (11.5,11.5);
\draw[out=60, in=185] (11.5,11.5) to (12.2,12.1);
\draw[out=5, in=205] (12.2,12.1) to (15,10.8);
\draw[out=25, in=245] (15,10.8) to (15.95,11.95);

%L(r)
\draw[out=30, in=190] (4.5,5.9) to (5.8,6.5);
\draw[out=10, in=200] (5.8,6.5) to (8,6.5);
\draw[out=20, in=220] (8,6.5) to (10.2,9.2);

\draw[out=40, in=170] (10.2,9.2) to (12,9.6);
\draw[out=-10, in=185] (12,9.6) to (13.4,9.4);
\draw[out=5, in=225] (13.4,9.4) to (14.8,10.3);
\draw[out=45, in=190] (14.8,10.3) to (15.95,10.95);

%\IM q_H(ir)
\draw[out=30, in=245] (4.5,6.05) to (5.7,7.5);
\draw[out=65, in=150] (5.7,7.5) to (7.5,8);
\draw[out=-30, in=240] (7.5,8) to (9.35,8.2);
\draw[out=60, in=221] (9.35,8.2) to (10.03,9.13);
\draw[out=41, in=237] (10.03,9.13) to (11.04,10.4);
\draw[out=57, in=180] (11.04,10.4) to (11.7,10.9);
\draw[out=0, in=140] (11.7,10.9) to (12.3,10.7);
\draw[out=-40, in=172] (12.3,10.7) to (13.52,9.75);
\draw[out=-8, in=218] (13.52,9.75) to (15,10.59);
\draw[out=38, in=241] (15,10.59) to (15.95,11.67);

%Legend
\node (A) at (4.5,5) {$\mr r(\xi_n)$};
\node (B) at (5.1,4.55) {$\hat r(\xi_n)$};
\node (C) at (5.85,5) {$\mr r(t_n)$};
\node (D) at (9.15,4.55) {$\hat r(t_n)$};
\node (E) at (9.9,5) {$\mr r(\xi_{n+1})$};

\node (F) at (10.5,4.55) {$\hat r(\xi_{n+1})$};
\node (G) at (11.25,5) {$\mr r(t_{n+1})$};
\node (H) at (13.7,4.55) {$\hat r(t_{n+1})$};
\node (I) at (14.8,5) {$\mr r(\xi_{n+2})$};
\node (J) at (15.4,4.55) {$\hat r(\xi_{n+2})$};

\draw[loosely dotted] (4.5,5.3) to (4.5,12.2);
\draw[loosely dotted] (5.1,4.85) to (5.1,12.2);
\draw[loosely dotted] (5.85,5.3) to (5.85,12.2);
\draw[loosely dotted] (9.15,4.85) to (9.15,12.2);
\draw[loosely dotted] (9.9,5.3) to (9.9,12.2);
\draw[loosely dotted] (10.5,4.85) to (10.5,12.2);
\draw[loosely dotted] (11.25,5.3) to (11.25,12.2);
\draw[loosely dotted] (13.7,4.85) to (13.7,12.2);
\draw[loosely dotted] (14.8,5.3) to (14.8,12.2);
\draw[loosely dotted] (15.55,4.85) to (15.55,12.2);

%tags
\node[scale=0.8] (J) at (7,6.28) {$L(r)$};
\node[scale=1.15] (K) at (7.1,7.15) {$r^{\frac{\log l}{\log (pl)}}$};
\node[scale=0.8] (L) at (7,8.48) {$\IM q_H(ir)$};
\node[scale=0.8] (M) at (7,9.18) {$A(r)$};

\end{tikzpicture}
A sketch of the behaviour of $q_{H_{p,l}}$
\end{center}

\noindent Formal details are given in the following theorem as well as in \Cref{R2}.

\begin{theorem}
\label{T2}
Let $p,l \in (0,1)$. For the Hamiltonian $H=H_{p,l}$ from \Cref{T6} and for all sufficiently large $n \in \bb N$, we have
\begin{equation}
\label{A8}
\mr r(\xi_n) < \hat r(\xi_n) <\mr r(t_n) < \hat r(t_n) < \mr r(\xi_{n+1}).
\end{equation}
On the intervals delimited by the terms in (\ref{A8}), the functions $L(r)$, $\IM q_H(ir)$, and $A(r)$ behave in the following way:
\begin{itemize}
\item[$(i)$] $\IM q_H(ir) \asymp A(r)$ uniformly for $r \in [\mr r(\xi_n),\hat r(\xi_n)]$, $n \in \bb N$.
\item[$(ii)$] $\IM q_H(ir) \asymp A(r)$ uniformly for $r \in [\hat r(\xi_n),\mr r(t_n)]$, $n \in \bb N$. \\[1ex]
Moreover, $L(\mr r(t_n)) \ll A(\mr r(t_n))$.
\item[$(iii)$] $L(r) \ll A(r)$ uniformly for $r \in [\mr r(t_n), \hat r(t_n)]$, $n \in \bb N$. \\[1ex]
In addition, $L(\mr r(t_n)) \ll \IM q_H(i \mr r(t_n)) \asymp A(\mr r(t_n))$ as well as $L(\hat r(t_n)) \asymp \IM q_H(i \hat r(t_n)) \ll A(\hat r(t_n))$.
\item[$(iv)$] $L(r) \asymp \IM q_H(ir) \ll A(r)$ uniformly for $r \in [\hat r(t_n),\mr r(\xi_{n+1})]$, $n \in \bb N$.
\end{itemize}
\end{theorem}

\noindent The proof of this theorem involves some (partly tedious) computations that are partly contained in the forthcoming lemma. \\
The symbol $\approx$ should mean equality up to an additive term that is bounded in $n$ and $t$.

\begin{lemma}
\label{R4}
For the Hamiltonian $H_{p,l}$, the following formulae hold.  
\newline

\begin{tabular}{|lcl|}
\hline$\log \mr r(t_n)$ & $\mkern-10mu\approx$ & $\mkern-10mu-n^2 \frac{\log(pl)}{2}-n \frac{\log \big(\frac lp \big)}{2}$ \begin{minipage}[c][9mm][t]{0.1mm}%
\end{minipage}\\[1ex]
\hline
$\log \mr r(\xi_n)$ & $\mkern-10mu\approx$ &  $\mkern-10mu-n^2 \frac{\log(pl)}{2}+n \frac{\log (pl)}{2}$ \begin{minipage}[c][9mm][t]{0.1mm}%
\end{minipage}\\
\hline
\end{tabular}
$\mkern-12mu$
\begin{tabular}{|lcl|}
\hline
$\log \hat r(t_n)$ & $\mkern-10mu\approx$ &  $\mkern-10mu-n^2 \frac{\log (pl)}{2}-n \frac{\log l}{2}$ 
\begin{minipage}[c][9mm][t]{0.1mm}%
\end{minipage}\\
\hline
$\log \hat r(\xi_n)$ & $\mkern-10mu\approx$ &  $\mkern-10mu-n^2 \frac{\log(pl)}{2}+n \frac{\log (pl)}{2}$ 
\begin{minipage}[c][9mm][t]{0.1mm}%
\end{minipage}\\
\hline
\end{tabular}

\begin{tabular}{|lcll|}
\hline
$\log \mr r(t)$ & $\mkern-10mu \approx$ & $\mkern-10mu-n^2 \frac{\log(pl)}{2}-n \frac{\log \big(\frac lp \big)}{2}+\log f(t)$, & $\mkern-10mu t \in [t_n,\xi_n]$. 
\begin{minipage}[c][9mm][t]{0.1mm}%
\end{minipage}\\
\hline
$\log \mr r(t)$ & $\mkern-10mu\approx$ & $\mkern-10mu-n^2 \frac{\log(pl)}{2}-n \frac{\log (pl)}{2}+\log (1-f(t))$, & $\mkern-10mu t \in [\xi_{n+1},t_n]$. 
\begin{minipage}[c][9mm][t]{0.1mm}%
\end{minipage}\\
\hline
$\log \hat r(t)$ & $\mkern-10mu\approx$ & $\mkern-10mu-n^2 \frac{\log(pl)}{2}-n \frac{\log \big(\frac lp \big)}{2}+\log f(t)-\frac {\log (1-f(t))}{2}$, & $\mkern-10mu t \in [t_n,\xi_n]$. 
\begin{minipage}[c][9mm][t]{0.1mm}%
\end{minipage}\\
\hline
$\log \hat r(t)$ & $\mkern-10mu\approx$ & $\mkern-10mu-n^2 \frac{\log(pl)}{2}-n \frac{\log (pl)}{2}+\frac {\log (1-f(t))}{2}$, & $\mkern-10mu t \in [\xi_{n+1},t_n]$. 
\begin{minipage}[c][9mm][t]{0.1mm}%
\end{minipage}\\
\hline
\end{tabular}
\newline
\end{lemma}

\begin{proof}
First we calculate
\begin{align}
&\log(\mr r(t_n))=-\frac 12 \log [(\omega_1 \omega_2)(t_n)] =\frac 12 \int_{t_n}^{t_1} (\alpha_1(s)+\alpha_2(s)) \DD s \nonumber\\
&= \sum_{k=1}^{n-1} \Big( \int_{t_{k+1}}^{\xi_{k+1}} \frac{-f'(s)}{f(s)} \DD s + \int_{\xi_{k+1}}^{t_k} \frac{f'(s)}{1-f(s)} \DD s \Big) \nonumber\\
&= \sum_{k=1}^{n-1} \Big(\log(1-p^{k+1})-(k+1)\log l + \log(1-l^{k+1}) -k\log p \Big) \nonumber\\
\label{A23}
&\approx -n^2 \frac{\log(pl)}{2}-n \frac{\log \big(\frac lp \big)}{2}.
\end{align}
This also leads to
\begin{align*}
\log \hat r(t_n) &=-\frac 12 \log (1-f(t_n)^2)+\log \mr r(t_n) \approx -\frac 12 \log (1-f(t_n))+\log \mr r(t_n) \\
&\approx -n^2 \frac{\log (pl)}{2}-n \frac{\log l}{2}.
\end{align*}
If $t \in [t_n,\xi_n]$, then
\begin{align*}
\log \mr r(t)=\log \mr r(t_n)-\int_{t_n}^t \frac{-f'(s)}{f(s)} \DD s \approx -n^2 \frac{\log(pl)}{2}-n \frac{\log \big(\frac lp \big)}{2}+\log f(t).
\end{align*}
If $t \in [\xi_{n+1},t_n]$, then
\begin{align*}
&\log \mr r(t)=\log \mr r(t_n)+\int_t^{t_n} \frac{f'(s)}{1-f(s)} \DD s \\
&\approx -n^2 \frac{\log(pl)}{2}-n \frac{\log (pl)}{2}+\log (1-f(t)).
\end{align*}
By adding $-\frac 12 \log (1-f(t)^2) \approx -\frac 12 \log (1-f(t))$, the analogous formula for $\hat r(t)$ follows. Lastly,
\begin{align*}
\log \mr r(\xi_n) &\approx -n^2 \frac{\log(pl)}{2}-n \frac{\log \big(\frac lp \big)}{2}+\log f(\xi_n) \\
&\approx -n^2 \frac{\log(pl)}{2}+n \frac{\log (pl)}{2}.
\end{align*}
and
\begin{align*}
\log \hat r(\xi_n) &=-\frac 12 \log (1-f(\xi_n)^2)+\log \mr r(\xi_n) \approx \log \mr r(\xi_n).
\end{align*}
\end{proof}

\begin{proof}[Proof of \Cref{T2}]
It follows from \Cref{R4} that $\hat r(\xi_n)<\mr r(t_n)$ and $\hat r(t_n) < \mr r(\xi_{n+1})$ for large enough $n$. The remaining two inequalities in (\ref{A8}) follow from the basic fact that $\mr r(t) < \hat r(t)$ for all $t \in (0,\infty)$. \\
We will now prove $(i)-(iv)$ in reverse order.
\\[1.7ex]
\underline{$(iv)$:} $\xi_{n+1} \leq \mr t(r) \leq t_n$ and $\xi_{n+1} \leq \hat t(r) \leq t_n$. By \Cref{R4},
\begin{align*}
&-n^2 \frac{\log(pl)}{2}-n \frac{\log (pl)}{2}+\frac 12 \log \big(1-f(\hat t(r))\big) \approx \log \hat r(\hat t(r))=\log r \\
&=\log \mr r(\mr t(r)) \approx -n^2 \frac{\log(pl)}{2}-n \frac{\log (pl)}{2}+\log \big(1-f(\mr t(r))\big).
\end{align*}
Hence,
\begin{align*}
\frac{\IM q_H(ir)}{A(r)} \asymp \sqrt{1-f \big(\hat t(r)\big)^2} \asymp 1-f \big(\mr t(r)\big)^2 =\frac{L(r)}{A(r)}. 
\end{align*}
In addition, 
\[
\frac{L \big( \mr r(\xi_{n+1}) \big)}{A \big( \mr r(\xi_{n+1}) \big)} \asymp 1-f(\xi_{n+1})^2 \asymp 1,
\]
while
\[
\frac{L \big(\hat r(t_n) \big)}{A \big(\hat r(t_n) \big)} \asymp 1-f \big(\mr t(\hat r(t_n)) \big)^2 \asymp \sqrt{1-f(t_n)^2} = p^{\frac n2} \ll 1.
\]
\\[1.7ex]
\underline{$(iii)$:} $\xi_{n+1} \leq \mr t(r) \leq t_n$ and $t_n \leq \hat t(r) \leq \xi_n$. Thus
\begin{align*}
&-n^2 \frac{\log(pl)}{2}-n \frac{\log \big(\frac lp \big)}{2}+\log f(\hat t(r))-\frac 12 \log \big(1-f(\hat t(r))\big) \approx \log \hat r(\hat t(r)) \\
&=\log \mr r(\mr t(r)) \approx -n^2 \frac{\log(pl)}{2}-n \frac{\log (pl)}{2}+\log \big(1-f(\mr t(r))\big).
\end{align*}
Consequently,
\[
\frac 12 \log \big(1-f(\hat t(r))\big) \approx n \log p + \log f(\hat t(r))-\log \big(1-f(\mr t(r))\big),
\]
which implies
\[
\sqrt{1-f(\hat t(r))} \asymp p^n \frac{f(\hat t(r))}{1-f(\mr t(r))}.
\]
Let us check that the term $f(\hat t(r))$ can be neglected. Using that $f(\mr t(r)) \leq 1-p^n$, we get
\[
\sqrt{1-f(\hat t(r))} \lesssim f(\hat t(r))
\]
which is only possible if $f(\hat t(r))$ stays away from $0$. As $f(\hat t(r)) <1$, this means that $f(\hat t(r)) \asymp 1$, leading to
\[
\frac{\IM q_H(ir)}{A(r)} \asymp \sqrt{1-f(\hat t(r))} \asymp \frac{p^n}{1-f(\mr t(r))}.
\]
Hence, $\IM q_H(i\mr r(t_n)) \asymp A(\mr r(t_n))$. Looking back at case $(iv)$, we know that $\IM q_H(i\hat r(t_n)) \asymp L(\hat r(t_n)) \ll A(\hat r(t_n))$. In particular, since $\frac{L(r)}{A(r)}=1-f(\mr t(r))^2$ is increasing for $r$ in $[\mr r(t_n),\hat r(t_n)]$, we have $L(r) \ll A(r)$ uniformly on this interval.\\[1.7ex]
\underline{$(ii)$:} $t_n \leq \mr t(r) \leq \xi_n)$ and $t_n \leq \hat t(r) \leq \xi_n$, leading to
\begin{align*}
&-n^2 \frac{\log(pl)}{2}-n \frac{\log \big(\frac lp \big)}{2}+\log f(\hat t(r))-\frac 12 \log \big(1-f(\hat t(r))\big) \approx \log \hat r(\hat t(r)) \\
&=\log \mr r(\mr t(r)) \approx -n^2 \frac{\log(pl)}{2}-n \frac{\log \big(\frac lp \big)}{2}+\log f(\mr t(r)).
\end{align*}
Hence
\begin{align*}
\sqrt{1-f(\hat t(r))} \asymp \frac{f(\hat t(r))}{f(\mr t(r))} > f(\hat t(r)).
\end{align*}
In particular, $1-f(\hat t(r))$ stays away from $0$, which means that
\[
\frac{\IM q_H(ir)}{A(r)} \asymp \sqrt{1-f(\hat t(r))} \asymp 1.
\]
In other words, $\IM q_H(ir) \asymp A(r)$ uniformly for $r \in [\hat r(\xi_n), \mr r(t_n)]$. As we already know, $L(\mr r(t_n)) \ll \IM q_H(i\mr r(t_n)) \asymp A(\mr r(t_n))$.\\[1.7ex]
\underline{$(i)$:} $t_n \leq \mr t(r) \leq \xi_n$ and $\xi_n \leq \hat t(r) \leq t_{n-1}$. In this case
\begin{align*}
&-n^2 \frac{\log(pl)}{2}+n \frac{\log (pl)}{2}+\frac 12 \log \big(1-f(\hat t(r))\big) \approx \log \hat r(\hat t(r))=\log \mr r(\mr t(r)) \\
&\approx -n^2 \frac{\log(pl)}{2}-n \frac{\log \big(\frac lp \big)}{2}+\log f(\mr t(r)).
\end{align*}
Taking into account that $f(\mr t(r)) \geq l^n$ by definition, it follows that
\[
\frac{\IM q_H(ir)}{A(r)} \asymp \sqrt{1-f(\hat t(r))} \asymp \frac{f(\mr t(r))}{l^n} \asymp 1.
\]
Therefore, $\IM q_H(ir) \asymp A(r)$ uniformly for $r \in [\mr r(\xi_n),\hat r(\xi_n)]$. At the left end of this interval, we even have $L(\mr r(\xi_n)) \asymp A(\mr r(\xi_n))$ by case $(iv)$.
\end{proof}

\noindent Before we state our next result, we note that by definition of $H_{p,l}$,
\begin{equation}
\label{A13}
\liminf_{t \to 0} \frac{\det \Omega(t)}{(\omega_1 \omega_2)(t)}=\liminf_{t \to 0} \big(1-f(t)^2 \big)=0.
\end{equation}
In view of (\ref{A9}), we have $\liminf_{r \to \infty} \frac{\IM q_{H_{p,l}}(ir)}{|q_{H_{p,l}}(ir)|}=0$ and hence $\IM q_{H_{p,l}}(ir) \not\asymp |q_{H_{p,l}}(ir)|$. \\
Nevertheless, the following lemma shows that $|q_{H_{p,l}}(ir)|$ grows faster than a power. Recalling \Cref{T7}, this means that $|q(ir)| \gtrsim r^{\delta}$ for $r \to \infty$ is not a sufficient condition for $\IM q(ir) \asymp |q(ir)|$ as $r \to \infty$. Instead, we see that $|q(ir)|$ being positively increasing really means that not only does $|q(ir)|$ grow sufficiently fast, but also without oscillating too much.

\begin{lemma}
\label{R2}
Let $\delta := \frac{\log l}{\log (pl)} \in (0,1)$. Then
\begin{itemize}
\item[$\rhd$] $|q_{H_{p,l}}(ir)| \gtrsim r^{\delta}$, $r \to \infty$,
\item[$\rhd$] $|q_{H_{p,l}}(i\mr r(\xi_n))| \asymp \mr r(\xi_n)^{\delta}$.
\end{itemize}
\end{lemma}
\begin{proof}
We start the proof with calculating, for $t \in [t_n,\xi_n]$,
\begin{align*}
&\log \sqrt{ \frac{\omega_1(t)}{\omega_2(t)} }=\frac 12 \log \Big(\frac{\omega_1(t)}{\omega_2(t)} \Big)=\sum_{k=1}^{n-1} \int_{t_{k+1}}^{\xi_{k+1}} \frac{-f'(s)}{f(s)} \DD s + \int_{t_n}^t \frac{f'(s)}{f(s)} \DD s\\
&=\sum_{k=1}^{n-1} \Big(\log (1-p^{k+1})-(k+1)\log l \Big)+\log f(t)-\log (1-p^n) \\
&\approx -(n^2+n)\frac{\log l}{2}+\log f(t).
\end{align*}
Now we use our formula for $\log \mr r(t)$:
\begin{align}
&\log \sqrt{\frac{\omega_1(t)}{\omega_2(t)} } \nonumber \\
&\approx \frac{\log l}{\log(pl)}\log \mr r(t)+\frac 12 \bigg(\frac{\log(l)\log(\frac lp)}{\log(pl)} -\log l \bigg)n + \bigg(1-\frac{\log l}{\log(pl)} \bigg)\log f(t)\nonumber\\
\label{A22}
&=\frac{\log l}{\log(pl)}\log \mr r(t)+\frac{\log p}{\log(pl)}\big(\log f(t)-n \log l\big), \quad t \in [t_n,\xi_n].
\end{align}
Since $f$ was assumed to be monotone decreasing on $[t_n,\xi_n]$, and $\log f(\xi_n)=n \log l$, 
\[
\log \sqrt{\frac{\omega_1(t)}{\omega_2(t)} } \gtrapprox \frac{\log l}{\log(pl)}\log \mr r(t)=\delta \log \mr r(t),
\]
where $\gtrapprox$ indicates that the inequality holds up to an additive term that is bounded in $n$ and $t$. Therefore
\[
|q_{H_{p,l}}(i \mr r(t))| \asymp \sqrt{\frac{\omega_1(t)}{\omega_2(t)}} \gtrsim \mr r(t)^{\delta}, \quad t \in [t_n,\xi_n].
\]
Observing that $\frac{\omega_1}{\omega_2}$ is constant on $[\xi_{n+1},t_n]$ (since $\alpha_1-\alpha_2=0$ there), we obtain this estimate also for $t \in [\xi_{n+1},t_n]$:
\[
|q_{H_{p,l}}(i\mr r(t))| \asymp \sqrt{\frac{\omega_1(t)}{\omega_2(t)}}=\sqrt{\frac{\omega_1(\xi_{n+1})}{\omega_2(\xi_{n+1})}} \gtrsim \mr r(\xi_{n+1})^{\delta} \geq \mr r(t)^{\delta}.
\]
Finally, setting $t=\xi_n$ in (\ref{A22}) yields $|q_{H_{p,l}}(i\mr r(\xi_n))| \asymp \mr r(\xi_n)^{\delta}$.
\end{proof}

\vspace{7pt}

\begin{example}
\label{R3}
Let $H$ be as in \Cref{T6}, but $f(\xi_n)=1-l^{n-1}$ instead, where $l > \sqrt{p}$. Similarly to \Cref{T2}, one can show that
\[
L(\hat r(t_n)) \asymp \IM q_H(i \hat r(t_n)) \ll A(\hat r(t_n)).
\]
However, for our new Hamiltonian, 
\[
\lim_{t \to 0} \frac{\det \Omega(t)}{(\omega_1 \omega_2)(t)} = \limsup_{t \to 0} \frac{\det \Omega(t)}{(\omega_1 \omega_2)(t)} = \limsup_{t \to 0} \big( 1-f(t)^2 \big) =  0
\]
as opposed to (\ref{A13}).

\end{example}

\section{Reformulation for Krein strings}
\label{S6}

Recall that a \textit{Krein string} is a pair $S[L,\mathfrak{m}]$ consisting of a number $L \in (0,\infty ]$ and a nonnegative Borel measure $\mathfrak{m}$ on $[0,L]$, such that $\mathfrak{m}([0,t])$ is finite for every $t \in [0,L)$, and $\mathfrak{m}(\{L\})=0$. To this pair we associate the equation
\begin{equation}
y_+'(x)+z\int_{[0,x]} y(t)\DD \mathfrak{m}(t)=0, \quad \quad x \in [0,L),
\end{equation}
where $y_+'$ denotes the right-hand derivative of $y$, and $z$ is a complex spectral parameter. \\
For each string, we can construct a function $q_S$ called the \textit{principal Titchmarsh-Weyl coefficient} of the string (\cite{langer.winkler:1998} following \cite{kac.krein:1968}). This function belongs to the Stieltjes class, i.e., it is analytic on $\bb C \setminus [0,\infty)$, its imaginary part is nonnegative on $\bb C_+$, and its values on $(-\infty ,0)$ are positive. The correspondence between Krein strings and functions of Stieltjes class is bijective, as was shown by M.G.Krein. \newline

\noindent \Cref{A41} below is the reformulation of \Cref{T1} for the Krein string case. 

\begin{theorem}
\label{A41}
Let $S[L,\mathfrak{m}]$ be a Krein string and set
\begin{equation}
\delta (t) := \bigg(\int_{[0,t)} \xi^2 \DD \mathfrak{m}(\xi) \bigg)\cdot \bigg(\int_{[0,t)} \DD \mathfrak{m}(\xi) \bigg)-\bigg(\int_{[0,t)} \xi \DD \mathfrak{m}(\xi) \bigg)^2.
\end{equation}
for $t \in [0,L)$. Let 
\[
\hat \tau (r):=\inf \big\{t>0 \, \big| \, \frac{1}{r^2} \leq \delta (t) \big\}, \quad \quad r \in (0,\infty).
\]
We set
\begin{align}
f(r) :=\mathfrak{m}([0,\hat \tau (r)))+ \mathfrak{m}(\{\hat \tau (r)\}) \frac{\frac{1}{r^2}-\delta (\hat \tau (r))}{\delta (\hat \tau (r)+)-\delta (\hat \tau (r))}
\end{align}
if $\delta$ is discontinuous at $\hat \tau(r)$, and $f(r):=\mathfrak{m}([0,\hat \tau (r)))$ otherwise. Then 
\begin{align}
\IM q_S(ir) \asymp \frac{1}{rf(r)}, \quad \quad r \in (0,\infty ),
\end{align}
with constants independent of the string.
\end{theorem}

\noindent Before proving \Cref{A41}, we need to introduce the concept of dual strings as well as a Hamiltonian associated to a string. Writing
\[
m(t) := \mathfrak{m}([0,t)), \quad \quad t \in [0,L)
\]
we can define the dual string $S[\hat L , \hat{\mathfrak{m}}]$ of $S[L,\mathfrak{m}]$ by setting
\[
\hat L :=\left\{ \begin{array}{ll}
m(L) & \text{if } L+m(L)=\infty ,\\
\infty & \text{else}
\end{array} \right.
\]
and
\[
\hat m (\xi):=\inf \{t >0 \, | \, \xi \leq m (t)\}.
\]
The function $\hat m$ is increasing and left-continuous and thus gives rise to a nonnegative Borel measure $\hat{\mathfrak{m}}$. \newline

\noindent The Hamiltonian defined by
\begin{equation}
\label{A39}
H(t) := \left\{ \begin{array}{ll} \begin{pmatrix}
\hat m(t)^2 & \hat m(t) \\
\hat m(t) & 1
\end{pmatrix} & \text{if } t \in [0,\hat L], \\[3ex]
\begin{pmatrix}
1 & 0 \\
0 & 0
\end{pmatrix} & \text{if } \hat L + \int_0^{\hat L} \hat m(t)^2 \DD t <\infty ,\,\, \hat L<t<\infty 
\end{array} \right.
\end{equation}
then satisfies $q_S=q_H$, see e.g. \cite{kaltenbaeck.winkler.woracek:2007}.

\begin{proof}[Proof of \Cref{A41}]
In view of \Cref{T1} and the fact that $q_S=q_H$ for the Hamiltonian $H$ defined in (\ref{A39}), our task is to express $\hat t_H(r)$ in terms of the string. If $\delta (\hat \tau (r)) = \frac{1}{r^2}$, this is easy because of \cite[Corollary 3.4]{kaltenbaeck.winkler.woracek:2007} giving
\begin{align*}
\det \Omega_H(m (\hat \tau (r))) = \delta (\hat \tau (r))=\frac{1}{r^2}
\end{align*}
and hence $\hat t_H(r)=m (\hat \tau (r))$. \\
Otherwise, we have $\delta (\hat \tau (r))<\frac{1}{r^2}$ and $\delta (\hat \tau (r)+) \geq \frac{1}{r^2}$. Using again \cite[Corollary 3.4]{kaltenbaeck.winkler.woracek:2007}, we have
\begin{equation}
\label{A40}
\det \Omega_H(m (\hat \tau (r)))= \delta (\hat \tau (r)) < \frac{1}{r^2}, \quad \quad \det \Omega_H(m (\hat \tau (r)+))= \delta (\hat \tau (r)+) \geq \frac{1}{r^2}
\end{equation}
which tells us that $\hat t_H(r) \in \big(m (\hat \tau (r)),m (\hat \tau (r)+) \big]$. By \cite[Lemma 3.1]{kaltenbaeck.winkler.woracek:2007}, $\hat m$ is constant on this interval. Therefore, for $t \in \big(m (\hat \tau (r)),m (\hat \tau (r)+) \big]$,
\begin{align*}
\det \Omega_H(t) &=\bigg(\int_0^{m(\hat \tau (r))} \hat m (x)^2 \DD x  +\big(t-m(\hat \tau (r))\big)\hat m (t)^2 \bigg)\cdot t \\
&-\bigg(\int_0^{m(\hat \tau (r))} \hat m (x) \DD x  +\big(t-m(\hat \tau (r))\big)\hat m (t) \bigg)^2 = c_1(r)t+c_2(r)
\end{align*}
for some constants $c_1(r),c_2(r)$. Using (\ref{A40}), this leads to
\[
\det \Omega_H(t)=\delta (\hat \tau (r)) + \frac{t-m(\hat \tau (r))}{m(\hat \tau (r)+)-m(\hat \tau (r))} \big(\delta (\hat \tau (r)+)-\delta (\hat \tau (r))\big).
\]
If we equate this to $\frac{1}{r^2}$, we find that
\[
\hat t_H(r)=m(\hat \tau (r))+ \big(m(\hat \tau (r)+)-m(\hat \tau (r)) \big) \frac{\frac{1}{r^2}-\delta (\hat \tau (r))}{\delta (\hat \tau (r)+)-\delta (\hat \tau (r))}=f(r).
\]
Now we have $\omega_{H;2}(t)=\int_0^t h_2(s) \DD s =t$, and \Cref{T1} now shows
\[
\IM q_S(ir)=\IM q_H(ir) \asymp \frac{1}{r\hat t_H(r)}=\frac{1}{rf(r)}.
\]
\end{proof}

\appendix
\appendixpage
\section{A construction method for Hamiltonians with prescribed angle of $\boldsymbol{q_H}$}
\label{APP}
\setcounter{lemma}{0}

Let $H$ be a Hamiltonian on $(a,b)$. Assume for simplicity that $\hat a=a$, i.e., $a$ is not the left endpoint of an $H$-indivisible interval. As discussed at the beginning of \Cref{S3}, the behavior of 
\[
\frac{\det \Omega (t)}{(\omega_1\omega_2)(t)}=1 - \frac{\omega_3(t)^2}{(\omega_1\omega_2)(t)} >0
\]
towards the left endpoint $a$ corresponds to the angle of $q_H(ir)$ for $r \to \infty$. It is thus desirable to be able to construct examples of Hamiltonians with prescribed $\frac{\omega_3(t)}{\sqrt{(\omega_1\omega_2)(t)}}$, which is what we did in \Cref{T6}. We give now a general version of this idea. \\
The following result is formulated for Hamiltonians in limit circle case, making the statement cleaner. When we made use of this construction method in \Cref{T6}, we obtained a Hamiltonian in limit point case by simply appending an infinitely long indivisible interval.

\begin{proposition}
\label{T3}
Let $f$ be locally absolutely continuous on $(a,b]$ and such that $f(t) \in (-1,1)$ for all $t \in (a,b]$. Then there is a Hamiltonian, in limit circle case at $b$, with the properties
\begin{itemize}
\item[(i)] $a=\hat a$ is not the left endpoint of an $H$-indivisible interval, and
\item[(ii)] $f(t)=\frac{\omega_3(t)}{\sqrt{(\omega_1\omega_2)(t)}}$ for all $t \in (a,b]$.
\end{itemize}
In addition, let
\[
\Delta(f):=\frac{2|f'|}{1-\sgn (f')f}
\]
which is in $L_{loc}^1((a,b])$.
Then all possible choices for $(\omega_1\omega_2)(t)$ are given by functions of the form
\begin{align*}
\exp \Big(c-\int_t^b g(s) \DD s \Big)
\end{align*}
 where $c \in \bb R$, $g \in L_{loc}^1((a,b]) \setminus L^1((a,b])$ with $g(t) \geq \Delta(f)(t)$ and $g(t)>0$ for $t \in (a,b]$ a.e.
\end{proposition}
\begin{proof}
If $H$ is given and such that $(i),(ii)$ hold, then clearly $f$ is locally absolutely continuous and takes values in $(-1,1)$. \\
Let $f$ be as in the statement. Then clearly $f' \in L_{loc}^1((a,b])$. Also, the denominator of $\Delta (f)$ is locally bounded below by a positive number, and hence $\Delta(f) \in L_{loc}^1((a,b])$. We check the conditions that $\omega_1(t),\omega_2(t)$ must satisfy in order that they, together with $\omega_3(t):=\sqrt{(\omega_1\omega_2)(t)}f(t)$, give rise to a Hamiltonian through $h_i(t):=\omega_i'(t)$, $i=1,2,3$. Clearly, $\omega_1,\omega_2$ have to be increasing, absolutely continuous on $[a,b]$ and satisfy $\omega_1(0)=\omega_2(0)=0$. Moreover, we want
\[
(h_1h_2)(t) \geq h_3(t)^2=\Big(\sqrt{(\omega_1\omega_2)(t)}f'(t)+\frac{h_1(t)\omega_2(t)+\omega_1(t)h_2(t)}{2 \sqrt{(\omega_1\omega_2)(t)}}f(t) \Big)^2.
\]
This is equivalent to
\[
\frac{(h_1h_2)(t)}{(\omega_1\omega_2)(t)} \geq \bigg(f'(t)+\frac 12 \Big(\frac{h_1(t)}{\omega_1(t)}+\frac{h_2(t)}{\omega_2(t)} \Big) f(t)\bigg)^2
\]
Setting 
\[
\alpha_i(t):=\frac{h_i(t)}{\omega_i(t)}, \quad i=1,2, \quad \quad g:=\alpha_1+\alpha_2 \in L_{loc}^1((a,b]),
\]
the inequality takes the form
\[
\alpha_i (g-\alpha_i) \geq \Big(f'+\frac 12 gf \Big)^2
\]
which is equivalent to 
\begin{align}
\label{A7}
\alpha_i \in \Bigg[\frac g2-\sqrt{\frac{g^2}{4}-\Big(f'+\frac 12 gf \Big)^2},\frac g2+\sqrt{\frac{g^2}{4}-\Big(f'+\frac 12 gf \Big)^2}\Bigg].
\end{align}
In particular, 
\begin{align*}
&\frac{g^2}{4}-\Big(f'+\frac 12gf \Big)^2 \geq 0 \, \Longleftrightarrow \, \frac g2 \geq \Big|f'+\frac 12gf \Big| \\
&\, \Longleftrightarrow g \geq \frac{2f'}{1-f} \text{ and } g \geq \frac{-2f'}{1+f} \\
&\, \Longleftrightarrow g \geq \frac{2|f'|}{1-\sgn(f')f} =\Delta(f).
\end{align*}
Since $\Delta(f) \in L_{loc}^1((a,b])$, we can find $g \in L_{loc}^1((a,b])$, $g(t) \geq \Delta (f)(t)$ a.e. on $(a,b]$, and additionally, $g \not\in L^1((a,b])$.
Choose measurable functions $\alpha_1,\alpha_2$ such that 
\begin{itemize}
\item[$\rhd$] $\alpha_1+\alpha_2=g$,
\item[$\rhd$] (\ref{A7}) holds for $\alpha_1$ at almost all $t \in (a,b)$ (and hence for $\alpha_2$), and
\item[$\rhd$] $\alpha_1,\alpha_2 \not\in L^1((a,b])$.
\end{itemize}
Note that $\alpha_1$ and $\alpha_2$ belong to $L_{loc}^1((a,b])$ since $g$ does, and that such a choice is possible because one can always take $\alpha_1=\alpha_2=\frac g2$. \\
From the construction it is clear that for a Hamiltonian $H$ with $\frac{d}{dt} [\log \omega_i(t)]=\alpha_i(t)$, $i=1,2$, there is $c \in \bb R$ such that
\[
\omega_i(t):=\exp \Big(c-\int_t^b \alpha_i(s) \DD s \Big), \quad \quad i=1,2.
\]
\end{proof}

\section{Calculations for \Cref{A24}}
\label{APPB}

Let $H$ be the Hamiltonian from \Cref{A24},
\begin{align*}
H(t)=
t^{\alpha -1}\left(\begin{matrix}
|\log t|^{\beta_1} & |\log t|^{\beta_3} \\
|\log t|^{\beta_3} & |\log t|^{\beta_2} \\
\end{matrix} \right), \quad \quad t \in (0,\infty),
\end{align*}
where $\alpha > 0$, $\beta_1, \beta_2 \in \bb R$ such that $\beta_1 \neq \beta_2$, and $\beta_3 := \frac{\beta_1 + \beta_2}{2}$. We carry out the calculations to justify the claimed asymptotics from the example. They were communicated by Matthias Langer. \\

 In order to calculate $\mr t(r)$ and $\hat t(r)$, two lemmas are needed.

\begin{lemma}
\label{approx_inv}
Let $f: (0,\varepsilon) \to (0,\infty)$ be increasing and $f(t) \sim ct^a |\log t|^b$ as $t \to 0$, for $a>0$, $c>0$, $b \in \bb R$. Then
\[
f^{-1}(s) \sim \Big(ca^{-b}s|\log s|^{-b} \Big)^{\frac 1a}, \quad s \to 0.
\]
\end{lemma}
\begin{proof}
We have
\begin{align}
\label{sim}
\lim_{t \to 0} \,\frac{f(t)}{t^a|\log t|^b}  = c.
\end{align}
Therefore,
\[
\lim_{t \to 0} \Big[\log f(t)-a\log t-b \log |\log t| \Big]= \log c
\]
and further
\[
\lim_{t \to 0} \Big[\frac{\log f(t)}{\log t} - a \Big]= \lim_{t \to 0} \Big[\frac{\log f(t)}{\log t} - a - b \frac{\log |\log t|}{\log t} \Big]= 0.
\]
In other words, $|\log t| \sim \frac 1a |\log f(t)|$. At the same time, by (\ref{sim}),
\[
t \sim \Big(c^{-1} f(t) |\log t|^{-b} \Big)^{\frac 1a} \sim \Big(c^{-1} f(t) [\frac 1a |\log f(t)|]^{-b} \Big)^{\frac 1a}
\]
which implies the assertion.
\end{proof}

\begin{lemma}
\label{int_asy}
Let $a>-1$ and $b \in \bb R$. Then
\begin{align*}
\int_0^t &s^a (-\log s)^b \DD s = \frac{1}{a+1}t^{a+1} (-\log t)^b \\
&\cdot \Big[1-\frac{b}{a+1}(-\log t)^{-1} + \frac{b(b-1)}{(a+1)^2}(-\log t)^2+\BigO \big((-\log t)^{-3} \big) \Big]
\end{align*}
\end{lemma}
\begin{proof}
\begin{align}
&\int_0^t s^a (-\log s)^b \DD s=\frac{1}{a+1}s^{a+1} (-\log s)^b \Big|_0^t - \frac{1}{a+1} \int_0^t s^a b (-\log s)^{b-1} \DD s \nonumber \\
\label{part_int_sa_logs_b}
&= \frac{1}{a+1}t^{a+1} (-\log t)^b - \frac{b}{a+1} \int_0^t s^a (-\log s)^{b-1} \DD s
\end{align}
Using (\ref{part_int_sa_logs_b}) two more times:
\begin{align*}
&= \frac{1}{a+1}t^{a+1} (-\log t)^b \\
&\quad- \frac{b}{a+1} \Big[\frac{1}{a+1}t^{a+1} (-\log t)^{b-1} - \frac{b-1}{a+1}\int_0^t s^a (-\log s)^{b-2} \DD s \Big] \\
&=\frac{1}{a+1}t^{a+1} (-\log t)^b \Big[1-\frac{b}{a+1}(-\log t)^{-1}+\frac{b(b-1)}{(a+1)^2}(-\log t)^{-2}\Big] \\
&\quad+c(a,b)\int_0^t s^a (-\log s)^{b-3} \DD s.
\end{align*}
The assertion follows using Karamata's Theorem \cite[Prop. 1.5.8 and 1.5.9a]{bingham.goldie.teugels:1989}.
\end{proof}

 We are now in position to determine $L(r)$, $\IM q_H(ir)$, and $A(r)$. 
\newline

\noindent \underline{Calculation of $\mr t(r)$:}
 By Karamata's Theorem we have
\[
\omega_i(t) \sim \frac{1}{\alpha}t^{\alpha} (-\log t)^{\beta_i}, \quad \quad i=1,2,3.
\]
Hence
\begin{align}
\label{m1m2}
(\omega_1\omega_2)(t) \sim \frac{1}{\alpha ^2}t^{2\alpha} (-\log t)^{\beta_1+\beta_2}.
\end{align}
Applying \Cref{approx_inv} yields
\[
(\omega_1\omega_2)^{-1}(s) \sim c \cdot s^{\frac{1}{2\alpha}} (-\log s)^{-\frac{\beta_3}{\alpha}}.
\]
We arrive at
\begin{align}
\label{tring}
\mr t(r)=(\omega_1\omega_2)^{-1}(r^{-2}) \sim c' \cdot r^{-\frac{1}{\alpha}} (\log r)^{-\frac{\beta_3}{\alpha}}.
\end{align}

\noindent \underline{Calculation of $A(r)$:}
\begin{align*}
A(r)&= \sqrt{\frac{\omega_1(\mr t(r))}{\omega_2(\mr t(r))}} \sim (-\log \mr t(r))^{\frac{\beta_1-\beta_2}{2}} \sim \Big[-\log \Big(r^{-\frac{1}{\alpha}} (\log r)^{-\frac{\beta_3}{\alpha}} \Big) \Big]^{\frac{\beta_1-\beta_2}{2}} \\
&\sim (\alpha \log r)^{\frac{\beta_1-\beta_2}{2}}.
\end{align*}

\noindent \underline{Calculation of $\hat t(r)$:} We use \Cref{int_asy} to calculate $\omega_i(t)$ with more precision:
\begin{align*}
\omega_i(t)=&\int_0^t s^{\alpha-1} (-\log s)^{\beta_i} \DD s = \frac{1}{\alpha}t^{\alpha} (-\log t)^{\beta_i} \\
&\cdot \Big[1-\frac{\beta_i}{\alpha}(-\log t)^{-1}+\frac{\beta_i (\beta_i -1)}{\alpha ^2}(-\log t)^{-2}+\BigO \big((-\log t)^{-3} \big) \Big]
\end{align*}
We get
\begin{align}
&\det \Omega (t)=\frac{1}{\alpha^2}t^{2\alpha} (-\log t)^{\beta_1+\beta_2} \nonumber\\
&\cdot \Big[1-\frac{\beta_1+\beta_2}{\alpha}(-\log t)^{-1}+\frac{\beta_1 (\beta_1-1)+\beta_2 (\beta_2-1) + \beta_1 \beta_2}{\alpha^2}(-\log t)^{-2} \nonumber\\
&-\Big(1-\frac{2\beta_3}{\alpha}(-\log t)^{-1}+\frac{2\beta_3 (\beta_3-1) + \beta_3^2}{\alpha^2}(-\log t)^{-2}+\BigO \big((-\log t)^{-3} \big) \Big)\Big] \nonumber\\
&=\frac{1}{\alpha^2}t^{2\alpha} (-\log t)^{\beta_1+\beta_2} \Big[\Big(\frac{\beta_1-\beta_2}{2\alpha} \Big)^2 (-\log t)^{-2} + \BigO \big((-\log t)^{-3} \big) \big)\Big] \nonumber\\
\label{detm}
&\sim c \cdot t^{2\alpha} (-\log t)^{2(\beta_3-1)}.
\end{align}
By \Cref{approx_inv},
\[
\big(\det \Omega \big)^{-1}(s) \sim c' \cdot s^{\frac{1}{2\alpha}} (-\log s)^{-\frac{\beta_3-1}{\alpha}}
\]
and further
\[
\hat t(r) = \big(\det \Omega \big)^{-1}(r^{-2}) \sim c'' \cdot r^{-\frac{1}{\alpha}} (\log r)^{-\frac{\beta_3-1}{\alpha}}.
\]

\noindent \underline{Calculation of $\IM q_H(ir)$:}
\begin{align*}
\IM q_H(ir)&\asymp \frac{1}{r\omega_2(\hat t(r))} \sim \frac{\alpha}{r} \big(\hat t(r) \big)^{-\alpha} \big(-\log \hat t(r) \big)^{-\beta_2} \\
&\sim c''' \cdot  \frac{\alpha}{r} r (\log r)^{\beta_3-1} (\log r)^{-\beta_2} = c''' \cdot (\log r)^{\frac{\beta_1-\beta_2}{2}-1}.
\end{align*}

\noindent \underline{Calculation of $L(r)$:} Using (\ref{detm}), (\ref{m1m2}) and (\ref{tring}), we have
\[
\frac{\det \Omega (\mr t(r))}{(\omega_1\omega_2)(\mr t(r))} \sim \Big(\frac{\beta_1-\beta_2}{2\alpha} \Big)^2 (\alpha \log r)^{-2}.
\]
Multiplying by $A(r)$, we obtain
\[
L(r) \sim c (\log r)^{\frac{\beta_1-\beta_2}{2}-2}.
\]

\subsection*{Acknowledgements.}
\noindent This work was supported by the Austrian Science Fund and the Russian Federation for Basic Research (grant number
I-4600). Furthermore, I would like to thank my supervisor Harald Woracek for his support and the expertise he provided.

%---------
%   FINISH
%---------

\printbibliography

\end{document}